\documentclass[12pt]{amsart}
\usepackage{amssymb}
\usepackage{amsfonts}
\usepackage{amsthm}
\usepackage{amsthm}
\usepackage{amscd}
\numberwithin{equation}{section}
\theoremstyle{plain}
\newtheorem{theorem}{Theorem}[section]

\newtheorem{proposition}[theorem]{Proposition}

\newtheorem{lemma}[theorem]{Lemma}
\theoremstyle{definition}
\newtheorem{definition}[theorem]{Definition}
\theoremstyle{remark}
\newtheorem*{remark}{Remark}

\begin{document}
%
%

\newcommand{\M}{\mathcal{M}_{g,N+1}^{(1)}}
\newcommand{\Teich}{\mathcal{T}_{g,N+1}^{(1)}}
\newcommand{\T}{\mathrm{T}}
\newcommand{\corr}{\bf}
\newcommand{\vac}{|0\rangle}
\newcommand{\Ga}{\Gamma}
\newcommand{\new}{\bf}
\newcommand{\define}{\def}
\newcommand{\redefine}{\def}
\newcommand{\Cal}[1]{\mathcal{#1}}
\renewcommand{\frak}[1]{\mathfrak{{#1}}}
\newcommand{\refE}[1]{(\ref{E:#1})}
\newcommand{\refS}[1]{Section~\ref{S:#1}}
\newcommand{\refSS}[1]{Section~\ref{SS:#1}}
\newcommand{\refT}[1]{Theorem~\ref{T:#1}}
\newcommand{\refO}[1]{Observation~\ref{O:#1}}
\newcommand{\refP}[1]{Proposition~\ref{P:#1}}
\newcommand{\refD}[1]{Definition~\ref{D:#1}}
\newcommand{\refC}[1]{Corollary~\ref{C:#1}}
\newcommand{\refL}[1]{Lemma~\ref{L:#1}}
\newcommand{\R}{\ensuremath{\mathbb{R}}}
\newcommand{\C}{\ensuremath{\mathbb{C}}}
\newcommand{\N}{\ensuremath{\mathbb{N}}}
\newcommand{\Q}{\ensuremath{\mathbb{Q}}}
\renewcommand{\P}{\ensuremath{\mathbb{P}}}
\newcommand{\Z}{\ensuremath{\mathbb{Z}}}
\newcommand{\kv}{{k^{\vee}}}
\renewcommand{\l}{\lambda}
\newcommand{\gb}{\overline{\mathfrak{g}}}
\newcommand{\g}{\mathfrak{g}}
\newcommand{\gh}{\widehat{\mathfrak{g}}}
\newcommand{\ghN}{\widehat{\mathfrak{g}_{(N)}}}
\newcommand{\gbN}{\overline{\mathfrak{g}_{(N)}}}
\newcommand{\tr}{\mathrm{tr}}
\newcommand{\sln}{\mathfrak{sl}}
\newcommand{\gl}{\mathfrak{gl}}
\newcommand{\Hwft}{\mathcal{H}_{F,\tau}}
\newcommand{\Hwftm}{\mathcal{H}_{F,\tau}^{(m)}}

\newcommand{\car}{{\mathfrak{h}}}    
\newcommand{\bor}{{\mathfrak{b}}}    
\newcommand{\nil}{{\mathfrak{n}}}    
\newcommand{\vp}{{\varphi}}
\newcommand{\bh}{\widehat{\mathfrak{b}}}  
\newcommand{\bb}{\overline{\mathfrak{b}}}  
\newcommand{\Vh}{\widehat{\mathcal V}}
\newcommand{\KZ}{Kniz\-hnik-Zamo\-lod\-chi\-kov}
\newcommand{\TUY}{Tsuchia, Ueno  and Yamada}
\newcommand{\KN} {Kri\-che\-ver-Novi\-kov}
\newcommand{\pN}{\ensuremath{(P_1,P_2,\ldots,P_N)}}
\newcommand{\xN}{\ensuremath{(\xi_1,\xi_2,\ldots,\xi_N)}}
\newcommand{\lN}{\ensuremath{(\lambda_1,\lambda_2,\ldots,\lambda_N)}}
\newcommand{\iN}{\ensuremath{1,\ldots, N}}
\newcommand{\iNf}{\ensuremath{1,\ldots, N,\infty}}

\newcommand{\MgN}{\mathcal{M}_{g,N}} 
\newcommand{\MgNeki}{\mathcal{M}_{g,N+1}^{(k,\infty)}} 
\newcommand{\MgNeei}{\mathcal{M}_{g,N+1}^{(1,\infty)}} 
\newcommand{\MgNekp}{\mathcal{M}_{g,N+1}^{(k,p)}} 
\newcommand{\MgNk}{\mathcal{M}_{g,N}^{(k)}} 
\newcommand{\MgNekpp}{\mathcal{M}_{g,N+1}^{(k,p')}} 
\newcommand{\MgNekkpp}{\mathcal{M}_{g,N+1}^{(k',p')}} 
\newcommand{\MgNezp}{\mathcal{M}_{g,N+1}^{(0,p)}} 
\newcommand{\MgNeep}{\mathcal{M}_{g,N+1}^{(1,p)}} 
\newcommand{\MgNeee}{\mathcal{M}_{g,N+1}^{(1,1)}} 
\newcommand{\MgNezz}{\mathcal{M}_{g,N+1}^{(0,0)}} 
\newcommand{\MgNi}{\mathcal{M}_{g,N}^{\infty}} 
\newcommand{\MgNe}{\mathcal{M}_{g,N+1}} 
\newcommand{\MgNep}{\mathcal{M}_{g,N+1}^{(1)}} 
\newcommand{\MgNp}{\mathcal{M}_{g,N}^{(1)}} 
\newcommand{\Mgep}{\mathcal{M}_{g,1}^{(p)}} 
\newcommand{\mpt}{\Sigma,P_1,P_2,\ldots, P_N,\Pif} 
\newcommand{\mpp}{\Sigma,P_1,P_2,\ldots, P_N} 
\newcommand{\At}{\widetilde{\Cal A}}
\newcommand{\Lt}{\widetilde{\Cal L}}
\newcommand{\Lc}{\overline{\Cal L}}
\newcommand{\Yt}{\widetilde{Y}}
\newcommand{\bt}{\tilde b}
\newcommand{\sinf}{{\widehat{\sigma}}_\infty}
\newcommand{\St}{\widetilde{S}}
\newcommand{\uni}{\mathcal{U}}
\newcommand{\can}{\mathcal{K}}
\newcommand{\Wh}{\widehat{W}}
\newcommand{\Wt}{\widetilde{W}}
\newcommand{\MegN}{\mathcal{M}_{g,N+1}^{(1)}} 

\newcommand{\Pif} {P_{\infty}}
\newcommand{\Pinf} {P_{\infty}}
\newcommand{\PN}{\ensuremath{\{P_1,P_2,\ldots,P_N\}}}
\newcommand{\PNi}{\ensuremath{\{P_1,P_2,\ldots,P_N,P_\infty\}}}
\newcommand{\Fln}[1][n]{F_{#1}^\lambda}
\newcommand{\tang}{\mathrm{T}}
\newcommand{\Kl}[1][\lambda]{\can^{#1}}
\newcommand{\A}{\mathcal{A}}
\newcommand{\U}{\mathcal{U}}
\newcommand{\V}{\mathcal{V}}
\renewcommand{\O}{\mathcal{O}}
\newcommand{\Ae}{\widehat{\mathcal{A}}}
\newcommand{\Ah}{\widehat{\mathcal{A}}}
\newcommand{\La}{\mathcal{L}}
\newcommand{\Le}{\widehat{\mathcal{L}}}
\newcommand{\Lh}{\widehat{\mathcal{L}}}
\newcommand{\eh}{\widehat{e}}
\newcommand{\Da}{\mathcal{D}}
\newcommand{\kndual}[2]{\langle #1,#2\rangle}
\newcommand{\cins}{\frac 1{2\pi\mathrm{i}}\int_{C_S}}
\newcommand{\cinsl}{\frac 1{24\pi\mathrm{i}}\int_{C_S}}
\newcommand{\cinc}[1]{\frac 1{2\pi\mathrm{i}}\int_{#1}}
\newcommand{\cintl}[1]{\frac 1{24\pi\mathrm{i}}\int_{#1 }}
\newcommand{\w}{\omega}
\newcommand{\ord}{\operatorname{ord}}
\newcommand{\res}{\operatorname{res}}
\newcommand{\nord}[1]{:\mkern-5mu{#1}\mkern-5mu:}
\newcommand{\Fn}[1][\lambda]{\mathcal{F}^{#1}}
\newcommand{\Fl}[1][\lambda]{\mathcal{F}^{#1}}
\renewcommand{\Re}{\mathrm{Re}}

\renewcommand{\k}{{k}}
\newcommand{\ce}{{c}}
\newcommand{\ka}{{k}}
\newcommand{\npos}{{p}}

\define\ldot{\hskip 1pt.\hskip 1pt}
\define\ifft{\qquad\text{if and only if}\qquad}
\define\a{\alpha}
\redefine\d{\delta}
\define\w{\omega}
\define\ep{\epsilon}
\redefine\b{\beta} \redefine\t{\tau} \redefine\i{{\,\mathrm{i}}\,}
\define\ga{\gamma}
\define\cint #1{\frac 1{2\pi\i}\int_{C_{#1}}}
\define\cintta{\frac 1{2\pi\i}\int_{C_{\tau}}}
\define\cintt{\frac 1{2\pi\i}\oint_{C}}
\define\cinttp{\frac 1{2\pi\i}\int_{C_{\tau'}}}
\define\cinto{\frac 1{2\pi\i}\int_{C_{0}}}
\define\cinttt{\frac 1{24\pi\i}\int_C}
\define\cintd{\frac 1{(2\pi \i)^2}\iint\limits_{C_{\tau}\,C_{\tau'}}}
\define\cintdr{\frac 1{(2\pi \i)^3}\int_{C_{\tau}}\int_{C_{\tau'}}
\int_{C_{\tau''}}}
\define\im{\operatorname{Im}}
\define\re{\operatorname{Re}}
\define\res{\operatorname{res}}
\redefine\deg{\operatornamewithlimits{deg}}
\define\ord{\operatorname{ord}}
\define\rank{\operatorname{rank}}
\define\fpz{\frac {d }{dz}}
\define\dzl{\,{dz}^\l}
\define\pfz#1{\frac {d#1}{dz}}

\define\K{\Cal K}
\define\U{\Cal U}
\redefine\O{\Cal O}
\define\He{\text{\rm H}^1}
\redefine\H{{\mathrm{H}}}
\define\Ho{\text{\rm H}^0}
\define\A{\Cal A}
\define\Do{\Cal D^{1}}
\define\Dh{\widehat{\mathcal{D}}^{1}}
\redefine\L{\Cal L} \redefine\D{\Cal D^{1}}
\define\KN {Kri\-che\-ver-Novi\-kov}
\define\Pif {{P_{\infty}}}
\define\Uif {{U_{\infty}}}
\define\Uifs {{U_{\infty}^*}}
\define\KM {Kac-Moody}
\define\Fln{\Cal F^\lambda_n}
\define\gb{\overline{\mathfrak{ g}}}
\define\G{\overline{\mathfrak{ g}}}
\define\Gb{\overline{\mathfrak{ g}}}
\redefine\g{\mathfrak{ g}}
\define\Gh{\widehat{\mathfrak{ g}}}
\define\gh{\widehat{\mathfrak{ g}}}
\define\Ah{\widehat{\Cal A}}
\define\Lh{\widehat{\Cal L}}
\define\Ugh{\Cal U(\Gh)}
\define\Xh{\hat X}
\define\Tld{...}
\define\iN{i=1,\ldots,N}
\define\iNi{i=1,\ldots,N,\infty}
\define\pN{p=1,\ldots,N}
\define\pNi{p=1,\ldots,N,\infty}
\define\de{\delta}

\define\kndual#1#2{\langle #1,#2\rangle}
\define \nord #1{:\mkern-5mu{#1}\mkern-5mu:}
\define \sinf{{\widehat{\sigma}}_\infty}
\define\Wt{\widetilde{W}}
\define\St{\widetilde{S}}
\define\Wn{W^{(1)}}
\define\Wtn{\widetilde{W}^{(1)}}
\define\btn{\tilde b^{(1)}}
\define\bt{\tilde b}
\define\bn{b^{(1)}}
%
\define\eps{\varepsilon}    
\define\doint{({\frac 1{2\pi\i}})^2\oint\limits _{C_0}
       \oint\limits _{C_0}}                            
\define\noint{ {\frac 1{2\pi\i}} \oint}   
\define \fh{{\frak h}}     
\define \fg{{\frak g}}     
\define \GKN{{\Cal G}}   
\define \gaff{{\hat\frak g}}   
\define\V{\Cal V}
\define \ms{{\Cal M}_{g,N}} 
\define \mse{{\Cal M}_{g,N+1}} 
\define \tOmega{\Tilde\Omega}
\define \tw{\Tilde\omega}
\define \hw{\hat\omega}
\define \s{\sigma}
\define \car{{\frak h}}    
\define \bor{{\frak b}}    
\define \nil{{\frak n}}    
\define \vp{{\varphi}}
\define\bh{\widehat{\frak b}}  
\define\bb{\overline{\frak b}}  
\define\Vh{\widehat V}
\define\KZ{Knizhnik-Zamolodchikov}
\define\ai{{\alpha(i)}}
\define\ak{{\alpha(k)}}
\define\aj{{\alpha(j)}}

\newcommand{\wf}[1]{\addtocounter{corr}{1}
            \footnote[\arabic{corr}]{#1}
            \marginpar{\footnotesize{** \arabic{corr}\ **}}
            }

\vspace*{-1cm}
\hspace*{\fill} math.AG/0312040
\vspace*{2cm}

\title[ Wess-Zumino-Witten-Novikov theory]
     {The Wess-Zumino-Witten-Novikov theory,
              Knizhnik-Zamolodchikov equations,
               and Krichever-Novikov algebras, II}
\author[M. Schlichenmaier]{Martin Schlichenmaier}
\address[Martin Schlichenmaier]{Laboratoire de Math\'ematique,
Universit\'e du Luxembourg,
162 A, Avenue de la Faiencerie,
L-1511 Luxembourg, Grand Duch\'e du Luxembourg}
\email{schlichenmaier@cu.lu}
\author[O.K. Sheinman]{Oleg K. Sheinman}
\thanks{O.K.Sheinman is supported by the
RFBR projects 02-01-00803, 02-01-22004 and  by the program
"Nonlinear Dynamics and Solitons" of the Russian Academy of
Science}
\address[ Oleg K. Sheinman]{Steklov Mathematical Institute, ul. Gubkina, 8,
Moscow 117966, GSP-1, Russia
and Independent University of Moscow,
Bolshoi Vlasievskii per. 11, Moscow 123298, Russia}
\email{oleg@sheinman.mccme.ru}

\begin{abstract}
This paper continues the same-named article, Part I \cite{rSSpt}.
We give a global operator approach to the WZWN theory for  compact
Riemann surfaces of an arbitrary genus $g$ with marked points.
Globality means here that we use Krichever-Novikov algebras of
gauge and conformal symmetries (i.e. algebras of global
symmetries) instead of loop and Virasoro algebras (which are local
in this context). The elements of this global approach are
described in Part I. In the present paper we give the construction
of conformal blocks and the projective flat connection on the
bundle constituted by them.
\end{abstract}
\subjclass{17B66, 17B67, 14H10, 14H15, 14H55,  30F30, 81R10, 81T40}
\keywords{
Wess-Zumino-Witten-Novikov theory, conformal blocks,
Knizhnik-Zamolodchikov equations, infinite-dimensional
Lie algebras, current algebras, gauge algebra, conformal
algebra, central extensions, highest weight representations,
Sugawara construction,  flat connection}
\date{31.11.2003}
\maketitle
\newpage
\tableofcontents
\section{Introduction}\label{S:intro}
In this article, we consider the following problem of
two-dimensional conformal field theory: induced by  gauge and
conformal symmetry  construct a vector bundle equipped with a
projective flat connection on the  moduli space of punctured
Riemann surfaces\footnote{See \refS{kzgen} for the precise
definition of what me mean by ``moduli space of punctured Riemann
surfaces''.}.

This problem originates in the well-known article of V.Knizhnik
and A.Za\-mo\-lod\-chi\-kov \cite{rKnZ}, where the case of genus
zero is considered. There the following remarkable system of
differential equations is obtained:
\begin{equation}\label{E:KZeqi}
  \bigg(k\frac{\partial}{\partial z_p}-
     \sum\limits_{r\ne p}\frac{t^a_pt^a_r}{z_p-z_r}\bigg)
    \Psi =0\ ,
                                    \qquad p=1,\ldots,N\ .
\end{equation}
Here $z_1,\ldots, z_N$ are arbitrary (generic) marked points on
the Riemann sphere. For $i=1,\ldots,N$,  representations $t_i$ of
a certain reductive Lie algebra $\g$ are given ($t_i^a$ being a
representation matrix for the $a^{th}$ generator of $\g$) and $k$
is a constant. A summation over $a$ is assumed in \refE{KZeqi}.

Nowadays, the equations are  known as the Knizhnik-Zamolodchikov
(KZ) equations. They can be interpreted as horizontality
conditions with respect to the \emph{Knizhnik-Zamolodchikov
connection}. From the point of view of physics, \refE{KZeqi} are
equations for the $N$-point correlation functions in
Wess-Zumino-Witten-Novikov models.

Further developments of the Knizhnik-Zamolodchikov ideas are
briefly outlined in \cite{rSSpt}. There more references can be
found. Until 1987 they  were inseparably linked with the
conception of gauge and conformal symmetries based on Kac-Moody
and Virasoro algebras. The higher genus generalizations in the
frame of this conception are initiated by D.~Bernard \cite{rBeRs},
the complete theory is given by A.~Tsuchiya, K.~Ueno, Y.~Yamada
\cite{rTUY}, see also \cite{EnFe}.

In \cite{rKNFa,rKNFb,rKNFc}, I.M.Krichever and S.P.Novikov defined
the basic objects of two-dimensional conformal field theory (like
the energy-momentum tensor) as global meromorphic objects on a
Riemann surface. For an abelian $\g$, and two punctures, they
pointed out another choice for the basic gauge and conformal
symmetries which are of a global nature and satisfy the
\emph{Krichever-Novikov algebras}.

In \cite{rSSS,rSSpt} we generalized the Krichever-Novikov results
to the non-abelian multi-point case and developed our \emph{global
operator approach} to the Wess-Zumino-Witten-Novikov model. In
\cite{rSSpt}, we basically formulated our approach including the
general form of the Knizhnik-Zamolodchikov connection for
arbitrary finite genus $g$ and the particular forms for lower
genera ($g=0$ or $g=1$).

We believe that the global operator approach simplifies the theory
and makes more transparent its geometry and relations. It enables
us  to  describe explicitly the Kuranishi tangent space of the
moduli space in terms of Krichever-Novikov basis elements, hence
to give explicitly  the equations of the generalized
Knizhnik-Zamolodchikov system. It is well-known that the higher
genus theories are related to the representations of the
fundamental group of the punctured Riemann surface ("twists" in
early CFT terminology \cite{FeWi}). In our approach,
representations of the fundamental group naturally arise as
parameters giving representations of gauge Krichever-Novikov
algebras. Thus, the global geometric Langlands correspondence
appears in the very beginning of the theory. One more intriguing
relation can  easily be formulated in the framework of our
approach, namely, the relation between Knizhnik-Zamolodchikov and
Hitchin systems \cite{rHit}. Both a Knizhnik-Zamolodchikov
operator and the corresponding quantized Hitchin integral
correspond to a Kuranishi tangent vector. On one hand, consider
the pull-back $e$ of such vector to the Krichever-Novikov vector
field algebra. On the other hand, consider the corresponding
derivative on the moduli space. In both cases subtract the
Sugawara operator corresponding to $e$. In the first case this
yields the quantized Hitchin integral \cite{ShN65}, and in the
second case the Knizhnik-Zamolodchikov operator.

However, in \cite{rSSpt}, we were not able to show that our
connection is well-defined on \emph{conformal blocks}, hence, we
also omitted the proof of its projective flatness. Filling up
these gaps is the goal of the present paper.

In \refS{kzkn} and \refS{kzrep} the necessary setup is revisited
from the point of view of the recent progress in the theory of
multi-point \KN\ algebras \cite{SchlCocycl, SchlMMJ},  and their
representations \cite{rShf,ShMMJ,ShUMN}. Certain results are
extended to be applicable to the situations considered here.

The \refS{kzgen} contains the main results, namely the
construction of the generalized Knizhnik-Zamolodchikov connection
on the conformal block bundle on an open dense part of the moduli
space, and the proof of the projective flatness of this
connection.

One of the authors (O.K.S.) was supported by the RFBR projects
02-01-00803, 02-01-22004 and  by the program "Nonlinear Dynamics
and Solitons" of the Russian Academy of Science. The authors thank
also the DFG for support in the framework of the
DFG-Forschergruppe ``Arithmetik'', located at the Universities of
Heidelberg and Mannheim.
\section{The algebras of Krichever-Novikov type}\label{S:kzkn}


For the general set-up developed in \cite{rSDiss},
\cite{rSLa,rSLb,rSLc} we refer to  \cite{rSSpt}. Let us introduce
here some notation.

Let $\Sigma$ be a compact Riemann surface of genus $g$, or in
terms of algebraic geometry, a smooth projective curve over $\C$
respectively. Let
$$
I=(P_1,\ldots,P_N),\quad  N\ge 1
$$
be a tuple of ordered, distinct points (``marked points''
``punctures'') on $\Sigma$ and $\Pif$ a distinguished marked point
on $\Sigma$  different from $P_i$ for every $i$. The points in $I$
are called the {\it in-points}, and the point $\Pif$ the {\it
out-point}.  Let $A=I\cup \{\Pif\}$ as a set. In \cite{rSLc},
\cite{rSDiss}, the general case where  there exists also
 a finite set of out-points is considered.
The results presented in this section are also valid in the more
general context.

\subsection{The Lie algebras $\A$, $\gb$, $\L$, $\D$ and $\D_\g$}
                                   \label{SS:algs}
${ }$

First let  $\A:=\A(\Sigma,I,\Pif)$ be the associative algebra of
meromorphic functions on $\Sigma$ which are regular except at the
points $P\in A$. Let  ${\g}$ be a complex finite-dimensional Lie
algebra. Then
\begin{equation}\label{E:curalg}
   {\gb}=\g\otimes_{\mathbb C}{\A}
\end{equation}
is called the {\it Krichever-Novikov current algebra}
\cite{rKNFa,rShea,rSha,rShns}. The Lie bracket on $\gb$ is given
by the relations
\begin{equation}\label{E:curr}
[x\otimes A,y\otimes B]=[x,y]\otimes AB.
\end{equation}
We will often omit  the symbol $\otimes$ in our notation.

Let $\L$ denote the Lie algebra  of meromorphic vector fields on
$\Sigma$ which are allowed to have poles only at the points $P\in
A$ \cite{rKNFa,rKNFb,rKNFc}.

For the Riemann sphere ($g=0$) with quasi-global coordinate $z$,
$I=\{0\}$ and $\Pif=\infty$, the algebra $\A$ is the algebra of
Laurent polynomials, the  current algebra $\gb$ is the loop
algebra, and the vector field algebra $\L$ is the Witt algebra.
Sometimes, we refer to this case as the {\it classical situation}
for short.

The algebra $\L$  operates on the elements of $\A$ by taking the
(Lie) derivative. This enables  us to define the Lie algebra $\Do$
of first order differential operators as the semi-direct sum of
$\A$ and $\L$. As vector space $\Do=\A\oplus \L$. The Lie
structure is defined by
\begin{equation}
[(g,e),(h,f)]:=(e\ldot h-f\ldot g,[e,f]),\quad g,h\in\A,\
e,f\in\L.
\end{equation}
Here $e\ldot f$ denotes  taking the Lie derivative. In local
coordinates $e_|=\tilde e\fpz$ and $e\ldot f_|=\tilde e\cdot
\pfz{f}$.

Consider at last the Lie algebra  $\D_\g$ of differential
operators associated to  $\gb$,
 (i.e. the {\it algebra of Krichever-Novikov
differential operators}). As a linear space $\D_\g=\gb\oplus\L$.
The Lie structure is given by the Lie structures on $\gb$, on
$\L$, and the additional definition
\begin{equation}\label{E:comm}
  [e,x\otimes A]:=-[x\otimes A,e]:=x\otimes (e\ldot A).
\end{equation}
In particular, for $\g={\frak{gl}}(1)$ one obtains, as a special
case, $\D_\g=\D$.
\subsection{Meromorphic forms of weight $\lambda$ and
\KN\ duality} ${ }$

Let $\K$ be the canonical line bundle. Its associated sheaf of
local sections is the sheaf of holomorphic differentials.
Following the common practice we will usually not distinguish
between a line bundle and its associated invertible sheaf of local
sections. For every $\l\in\Z$ we consider the bundle $\
\K^\l:=\K^{\otimes \l}$. Here we follow  the usual convention:
$\K^0=\Cal O$ is the trivial bundle,  and $\K^{-1}=\K^*$ is the
holomorphic tangent line bundle. Indeed, after fixing a theta
characteristics, i.e. a bundle  $S$ with $S^{\otimes 2}=\K$, it is
possible to consider $\l\in \frac {1}{2}\Z$. Denote by $\Fl$ the
(infinite-dimensional) vector space of global meromorphic sections
of $\K^\l$ which are holomorphic on $\Sigma\setminus A$. The
elements of $\Fl$  are called (meromorphic) forms or tensors of
weight $\lambda$.

Special cases, which are of particular interest, are the functions
($\l=0$), the vector fields ($\l=-1$), the  1-forms ($\l=1$), and
the quadratic differentials ($\l=2$). The space of functions is
already denoted by $\A$, and the space of vector fields by $\L$.

By  multiplying  sections with functions we again obtain sections.
In this way the $\Fl$ become $\A$-modules. By taking the Lie
derivative of the forms with respect to the vector fields the
vector spaces  $\Fl$ become $\L$-modules.
 In local coordinates the Lie derivative is
given as
\begin{equation}\label{E:Lder}
(e\ldot g)_|:=(\tilde e(z)\fpz)\ldot (\tilde g(z)\dzl):=
\left(\tilde e(z)\pfz {\tilde g}(z)+\l\, \tilde g(z) \pfz {\tilde
e}(z)\right)\dzl \ .
\end{equation}
The vector spaces $\Fl$ become $\Do$-modules by the canonical
definition $(g+e)\ldot v=g\cdot v+e\ldot v$. Here $g\in\A$,
$e\in\L$ and $v \in\Fl$. By universal constructions algebras of
differential operators of arbitrary degree can be considered
\cite{rSDiss,rSCt}.

Let $\rho$ be a meromorphic differential which is holomorphic on
$\Sigma\setminus A$ with exact pole order $1$ at the points in $A$
and given positive residues at $I$ and given negative residues at
$\Pif$ (of course obeying the restriction $\sum_{P\in
I}\res_P(\rho)+ \res_{\Pif}(\rho)=0$) and purely imaginary
periods. There exists exactly one such $\rho$ (see
\cite[p.116]{rSRS}). For $R\in \Sigma\setminus A$ a fixed point,
the function $u(P)=\Re\int_R^P\rho$ is a well-defined harmonic
function. The family of level lines $ C_\tau:=\{p\in M\mid
u(P)=\tau\},\ \tau\in\R $ defines a fibration of $\Sigma\setminus
A$. Each $C_\tau$ separates the points in $I$ from the point
$\Pif$. For $\tau\ll 0$ ($\tau\gg 0$) each level line $C_\tau$ is
a disjoint union of deformed circles $C_i$ around the points
$P_i$, $i=1,\ldots,N$ (of a deformed circle $C_\infty$  around the
point $\Pif$). We will call any such level line or any cycle
{homologous} to such a level line a separating cycle $C_S$.
\begin{definition}\label{D:knpair}
The {\it Krichever-Novikov pairing} ({\it KN pairing}) is the
pairing between $\Fl$ and $\Fl[1-\l]$ given by
\begin{equation}\label{E:knpair}
\begin{gathered}
\Fl\times\Fn[1-\l]\ \to\ \C,
\\
\kndual {f}{g}:=\cins f\otimes g =\sum_{P\in I}\res_{P}(f\otimes
g)= -\res_{\Pif}(f\otimes g),
\end{gathered}
\end{equation}
where $C_S$ is any separating cycle.
\end{definition}
The last equality follows from the residue theorem. Note that in
\refE{knpair} the integral does not depend on the separating cycle
chosen. From the construction of special dual basis elements in
the next subsection it follows that the KN pairing is
non-degenerate.

\subsection{\KN\ bases}\label{SS:knb}
${ }$

Krichever and Novikov introduced special bases ({\it \KN\ bases})
for the  vector spaces of meromorphic  tensors on Riemann surfaces
with two marked points. For $g=0$ the Krichever-Novikov bases
coincide with the Laurent bases. The multi-point generalization of
these bases is given by one of the authors in \cite{rSLc,rSDiss}
(see also Sadov \cite{rSad} for some results in similar
direction).  We define here the Krichever-Novikov type bases for
tensors  of arbitrary weight $\lambda$ on Riemann surfaces with
$N$ marked points as introduced in \cite{rSLc,rSDiss}.

For fixed $\l$ and for  every $n\in\Z$, and $p=1,\ldots,N$ a
certain element $f_{n,p}^\l\in\Fl$ is exhibited. The basis
elements are chosen in such a way that they fulfill the duality
relation
\begin{equation}\label{E:edu}
\kndual {f_{n,p}^\l} {f_{m,r}^{1-\l}}= \de_{-n}^{m}\cdot
\de_{p}^{r}
\end{equation}
with respect to the KN pairing \refE{knpair}. In particular, this
implies that the KN pairing is non-degenerate. Additionally, the
elements fulfill
\begin{equation}\label{E:ordfn}
\ord_{P_i}(f_{n,p}^\l)=(n+1-\l)-\d_i^p,\quad i=1,\ldots,N .
\end{equation}
The recipe for choosing the order at the point $\Pif$ is such that
up to a scalar multiplication there is a unique such element which
also fulfills \refE{edu}. After choosing local coordinates $z_p$
at the points $P_p$ the scalar can be fixed by requiring
\begin{equation}
{f_{n,p}^\l}_|(z_p)=z_p^{n-\l}(1+O(z_p))\left(dz_p\right)^\l,
\quad p=1,\ldots,N.
\end{equation}
To give an impression about the other requirement let us consider
the case   $g\ge 2$, $\l\ne 0,1$ and a generic choice for the
points in $A$  (or $g=0$ without any restriction). Then we require
\begin{equation}\label{E:ordi}
\ord_{\Pif}(f_{n,p}^\l)=-N\cdot(n+1-\l) +(2\l-1)(g-1)\ .
\end{equation}
By Riemann-Roch type arguments, it is shown in \cite{rSLa} that
there exists only one such element. For the necessary modification
for other cases see \cite{rSSpt}, \cite{rSLc,rSDiss}.

For the basis elements $f_{n,p}^\l$, explicit descriptions in
terms of rational functions (for $g=0$), the Weierstra\ss\
$\sigma$-function (for $g=1$), and prime forms and theta functions
(for $g\ge1$) are given in \cite{rSLb}. For $g=0$ and $g=1$, such
a description  can be found also in \cite[\S\S 2,7]{rSSpt}. For a
description using Weierstra\ss\ $\wp$-function, see \cite{rRDS},
\cite{rSDeg}. The existence of such  a description is necessary in
our context because we want to consider the above algebras  and
modules over the configuration space, respectively, the moduli
space of curves with marked points. In particular, one observes
from the explicit representation that the elements vary
``analytically'' when the complex structure of the Riemann surface
is deformed.

For the following special cases we introduce the notation:
\begin{equation}\label{E:econc}
A_{n,p}:=f_{n,p}^0,\quad e_{n,p}:=f_{n,p}^{-1},\quad
\w^{n,p}:=f_{-n,p}^1,\quad \Omega^{n,p}:=f_{-n,p}^2 \ .
\end{equation}

For $g=0$ and $N=1$ the basis elements constructed coincide with
the standard generators of the Witt and loop algebras,
respectively. For $g\ge 1$ and $N=1$ these elements coincide up to
an index shift with those given by Krichever and Novikov
\cite{rKNFa,rKNFb,rKNFc}.

\subsection{Almost graded structure, triangular
  decompositions}\label{SS:ags}
${ }$

For $g=0$ and $N=1$ the Lie algebras introduced in \refSS{algs}
are graded. A grading is a necessary tool for developing their
structure theory and the theory of their highest weight
representations. For the higher genus case (and for the
multi-point situation for $g=0$) there is no grading. It is a
fundamental observation due to Krichever and Novikov
\cite{rKNFa,rKNFb, rKNFc} that a weaker concept, an almost
grading, is sufficient to develop a suitable structure and
representation theory in this more general context.

An (associative or Lie) algebra is called {\it almost-graded} if
it admits a direct decomposition as a vector space $\
\V=\bigoplus_{n\in\Z} \V_n\ $, with (1) $\ \dim \V_n<\infty\ $ and
(2) there are constants $R$ and  $S$ such that
\begin{equation}\label{E:eaga}
 \V_n\cdot \V_m\quad \subseteq \bigoplus_{h=n+m-R}^{n+m+S} \V_h,
 \qquad\forall n,m\in\Z\ .
\end{equation}
The elements of $\V_n$ are called {\it homogeneous  elements of
degree $n$}. Let $\V=\oplus_{n\in\Z} \V_n$ be an almost-graded
algebra and
 $M$ an $\V$-module.
The module $M$ is called an {\it almost-graded} $\V$-module if it
admits a direct decomposition as a vector space $\
M=\bigoplus_{m\in\Z} M_m\ $, with (1) $\ \dim M_m<\infty\ $ and
(2) there are constants $T$ and  $U$ such that
\begin{equation}\label{E:eagm}
 V_n\ldot M_m\quad \subseteq \bigoplus_{h=n+m-T}^{n+m+U} M_h,
 \qquad\forall n,m\in\Z\ .
\end{equation}
The elements of $M_n$ are called {\it homogeneous  elements of
degree $n$}.

In our case the homogeneous subspaces $\Fln$ are defined as the
subspace of $\Fl$ generated by the elements $f_{n,p}^\l$ for
$p=1,\ldots,N$. Then $ \Fl=\bigoplus_{n\in\Z}\Fln$.
\begin{proposition}\label{P:almgrad}\cite{rSLc,rSDiss}
With respect to the introduced degree, the vector field algebra
$\L$, the function algebra $\A$, and the differential operator
algebra $\D$  are almost-graded and the $\Fl$ are  almost-graded
modules over them.
\end{proposition}
\noindent

The algebra  $\A$ can be decomposed (as  vector space) as follows:
\begin{equation}\label{E:edeca}
\begin{gathered}
\A=\A_+\oplus\A_{(0)}\oplus\A_-,
\\
\A_+:=\langle A_{n,p}\mid n\ge 1,\pN\rangle,\quad \A_-:=\langle
A_{n,p}\mid n\le -K-1,\pN\rangle\ ,\quad
\\
\A_{(0)}:=\langle A_{n,p}\mid -K\le  n\le 0, \pN\rangle \ ,
\end{gathered}
\end{equation}
and the Lie algebra $\L$ as follows:
\begin{equation}\label{E:edecv}
\begin{gathered}
\L=\L_+\oplus\L_{(0)}\oplus\L_-,
\\
\L_+:=\langle e_{n,p}\mid n\ge 1,\pN\rangle,\quad \L_-:=\langle
e_{n,p}\mid n\le -L-1,\pN\rangle\ ,\quad
\\
\L_{(0)}:=\langle e_{n,p}\mid -L\le  n\le 0, \pN\rangle \ .
\end{gathered}
\end{equation}
We call \refE{edeca}, \refE{edecv} the {\it triangular
decompositions}. In a similar way we obtain a triangular
decomposition of $\D$.

Due to the almost-grading the subspaces $\A_\pm$ and $\L_\pm$ are
subalgebras but the subspaces $\A_{(0)}$, and $\L_{(0)}$ in
general are not. We use the term {\it critical strip} for  them.

Note that $\A_+$, resp. $\L_+$ can be described as the algebra of
functions (vector fields) having a zero of at least  order one
(two) at the points $P_i,\iN$. These algebras can be enlarged by
adding all elements which are regular at all $P_i$'s. This can be
achieved by moving the set of basis elements $\ \{A_{0,p},\pN\}\
$, (resp. $\ \{e_{0,p},e_{-1,p},\iN\}\ $ from the critical strip
to these algebras. We denote the enlarged algebras by $\A_+^*$,
resp. by $\L_+^*$.

On the other hand $\A_-$ and $\L_-$ could also be enlarged such
that they contain all elements  which are regular at $\Pif$. This
is explained in detail in \cite{rSSpt}. We obtain  $\A_-^*$ and
$\L_-^*$ respectively. In the same way for every $p\in\N_0$ let
$\L_-^{(p)}$ be the subalgebra of vector fields vanishing of order
$\ge p+1$ at the point $\Pif$, and $\A_-^{(p)}$  the subalgebra of
functions respectively vanishing of order $\ge p$ at the point
$\Pif$. We obtain a decomposition
\begin{equation}
\L=\L_+\oplus\L_{(0)}^{(p)}\oplus\L_-^{(p)},\ \text{ for $p\ge
0$}, \quad \text{and}\quad
\A=\A_+\oplus\A_{(0)}^{(p)}\oplus\A_-^{(p)},\ \text{ for $p\ge
1$},
\end{equation}
with ``critical strips'' $\L_{(0)}^{(p)}$ and $\A_{(0)}^{(p)}$,
which are only subspaces. Of particular interest to us is
$\L_{(0)}^{(1)}$ which we call reduced critical strip. For $g\ge
2$ its dimension is
\begin{equation}\label{E:dimcs}
\dim \L_{(0)}^{(1)}=N+N+(3g-3)+1+1=2N+3g-1\ .
\end{equation}
The first two terms correspond to $\L_0$ and $\L_{-1}$. The
intermediate term comes from the vector fields in the basis which
have poles at the $P_i,\iN$ and $\Pif$. The $1+1$ corresponds to
the  vector fields in the basis with  exact order zero (one) at
$\Pif$.

The almost-grading can easily be extended to the {\it higher genus
current algebra} $\gb$ by setting  $\deg(x\otimes A_{n,p}):=n$. We
obtain a triangular decomposition as above
\begin{equation}\label{E:espcur}
\G=\G_+\oplus\G_{(0)}\oplus\G_-,\quad\text{with}\quad
\G_\beta=\g\otimes \A_\beta,\quad \beta\in\{-,(0),+\}\ ,
\end{equation}
In particular, $\G_\pm$ are subalgebras. The corresponding is true
for the enlarged subalgebras. Among them, $
\gb^r:=\gb_-^{(1)}=\g\otimes \A_-^{(1)} $ is of special
importance. It is called the regular subalgebra.

The finite-dimensional Lie algebra $\g$ can  naturally be
considered as subalgebra of $\Gb$. It lies in the subspace $\G_0$.
To see this we use $1=\sum_{p=1}^N A_{0,p}$, see \cite[Lemma
2.6]{rSSpt}.


\subsection{Central extensions and 2-cohomologies}\label{SS:cocyc}
${ }$

Let $\V$ be a Lie algebra and $\ga$  a Lie algebra 2-cocycle on
$\V$, i.e. $\ga$ is an antisymmetric bilinear form obeying
\begin{equation}\label{E:cocycle}
\ga([f,g],h)+\ga([g,h],f)+\ga([h,f],g)=0,\quad\forall f,g,h\in\V.
\end{equation}
On $\widehat\V=\C\oplus \V$ a Lie algebra structure can be defined
by  (with the notation $\widehat{f}:=(0,f)$ and $t:=(1,0)$)
\begin{equation}\label{E:centext}
[\widehat{f},\widehat{g}]:=\widehat{[f,g]}+\ga(f,g)\cdot t,\quad
[t,\widehat\V]=0.
\end{equation}
The element $t$ is a central element. Up to equivalence central
extensions are classified by the elements of  $\H^2(\V,\C)$, the
second Lie algebra cohomology space with values in the trivial
module $\C$. In particular, two cocycles $\ga_1,\ga_2$ define
equivalent central extensions if and only if there exist a linear
form $\phi$ on $\V$ such that
\begin{equation}
\ga_1(f,g)=\ga_2(f,g)+\phi([f,g]).
\end{equation}

\begin{definition}\label{D:local}
Let $\V=\bigoplus_{n\in\Z} \V_n$ be an almost-graded Lie algebra.
A cocycle $\ga$ for $\V$ is called local (with respect to the
almost-grading) if there exist $M_1,M_2\in\Z$ with
\begin{equation}\label{E:local}
\forall n,m\in\Z:\quad \ga(\V_n,\V_m)\ne 0\implies M_2\le n+m\le
M_1.
\end{equation}
\end{definition}
By defining $\deg(t):=0$ the central extension $\widehat\V$ is
almost-graded if and only if it is given by a local cocycle $\ga$.
In this case we call $\widehat\V$ an {\it almost-graded central
extension} or a
 {\it local central extension}.

In the following we consider cocycles of geometric origin for the
algebra introduced above.  First we deal with  $\A$, $\L$ and
$\Do$.  A thorough treatment for them is given in
\cite{SchlCocycl}.  The proofs of the following statements and
more details  can be found there.

For the abelian Lie algebra $\A$ any antisymmetric bilinear form
will be a 2-cocycle. Let $C$ be any (not necessarily connected)
differentiable cycle in $\Sigma\setminus A$ then
\begin{equation}\label{E:fung}
\ga^{(f)}_C:\A\times \A\to\C,\quad \gamma_C^{(f)}(g,h):=\cinc{C}
gdh
\end{equation}
is antisymmetric, hence a cocycle. Note that replacing $C$ by any
homologous (differentiable) cycle one obtains the same cocycle.
The above cocycle is \emph{$\L$-invariant}, i.e.
\begin{equation}\label{E:linv}
\ga^{(f)}_C(e\ldot g,h)=\ga^{(f)}_C(e\ldot h,g), \quad \forall
e\in\L,\ \forall g,h\in\A.
\end{equation}

For the vector field algebra $\L$ we want, following \cite{rKNFa},
to generalize the standard Virasoro-Gelfand-Fuks cocycle to higher
genus. To this end, we have first to choose a projective
connection. It will allow us to add a counter term to the
integrand to  obtain a well-defined 1-differential. Let $R$ be a
global holomorphic projective connection, see e.g. \cite{rSSpt}
for the definition. For every cycle $C$ (and every $R$) a cocycle
is given by

\begin{equation}\label{E:ecva}
\ga^{(v)}_{C,R}(e,f) :=\cinttt\left(\frac 12(\tilde e'''\tilde
f-\tilde e\tilde f''') -R\cdot(\tilde e'\tilde f-\tilde e\tilde
f')\right)dz\ .
\end{equation}
Here $e_|=\tilde e\frac {d}{dz}$ and $f_|=\tilde f\frac {d}{dz}$
with local meromorphic functions $\tilde e$ and $\tilde f$. A
different choice of the projective connection (even if we allow
meromorphic projective connections with poles only at the points
in  $A$) yields a cohomologous cocycle, hence an equivalent
central extension.

These two types of cocycles can be extended to cocycles {of} $\Do$
by setting {them} to be zero if one of the entries is from the
complementary space. For the vector field cocycles this is clear,
for the function algebra cocycles the $\L$-invariance \refE{linv}
is crucial. But there are other independent types of cocycles
which mix functions with vector fields. To define them we first
have to fix an affine connection $T$ which is holomorphic outside
$A$ and has at most a pole of order one at $\Pif$. For the
definition and existence of an affine connection, see
\cite{rSDiss}, \cite{ShMMJ}, \cite{SchlCocycl}. Now
\begin{equation}\label{E:mixg}
\ga_{C,T}^{(m)}(e,g):=-\ga_{C,T}^{(m)}(g,e):= \cinc
{C}\left(\tilde e\cdot g''+T\cdot (\tilde e\cdot
 g')\right)dz
\end{equation}
is a 2-cocycle. Again, the cohomology class does not depend on the
chosen affine connection.

Next we consider cocycles obtained by integrating over a
separating cycle $C_S$. Instead of $\ga_{C_S}$ we will use
$\ga_S$. Clearly, these cocycles can be expressed via residues at
the points in $I$ or equivalently at the point $\Pif$.
\begin{proposition}\label{P:local}
\cite{rSDiss} The above cocycles if integrated over a separating
cycle $C_S$ are local. In each case the upper bound is given by
zero.
\end{proposition}
If we replace $R$ or $T$ by other meromorphic connections which
have poles only at $A$, the cocycles still will be local. The
upper and lower  bounds might change. If the poles are of at most
order two for the projective connection, or order one for the
affine connection at the points in $I$, then the upper bounds will
remain  zero. Note that such a change of the connection can always
be given by adding elements from $\mathcal{F}^2$ or
$\mathcal{F}^1$ to $R$ and $T$ respectively.

As explained above we obtain almost-graded extensions
 $\Ah$, $\Lh$ and $\Dh$ via these cocycles.
By the vanishing of the cocycles (at least if $R$ is holomorphic)
on the subalgebras $\A_{\pm}$ and $\L_{\pm}$ the subalgebras can
be identified in a natural way with the subalgebras $\Ah_{\pm}$
and $\Lh_{\pm}$ of $\Ah$, resp. $\Lh$.

One of the main results of \cite{SchlCocycl} is
\begin{theorem}\label{T:unique}
(a) Every local cocycle of $\A$ which is $\L$-invariant is a
multiple of the cocycle $\ga_{S}^{(f)}$. The cocycle
$\ga_{S}^{(f)}$ is cohomologically non-trivial.
\newline
(b) Every local cocycle of $\L$ is cohomologous to a scalar
multiple of $\ga^{(v)}_{S,R}$. The cocycle  $\ga^{(v)}_{S,R}$
defines a non-trivial cohomology class. For every cohomologically
non-trivial local cocycle a meromorphic projective connection $R'$
which is holomorphic outside of $A$ can be chosen such that the
cocycle is equal to a  scalar multiple of
 $\ga^{(v)}_{S,R'}$.
\newline
(c) Every local cocycle for $\Do$ is a linear combination of the
above introduced cocycles  $\ga_{S}^{(f)}$,  $\ga^{(v)}_{S,R}$ and
 $\ga_{S,T}^{(m)}$ up to coboundary, i.e.
\begin{equation}
\gamma=r_1 \ga_{S}^{(f)}+ r_2 \ga_{S,T}^{(m)}+
r_3\ga^{(v)}_{S,R}+\text{coboundary},\quad r_1,r_2,r_3\in\C.
\end{equation}
The 3 basic cocycles are linearly independent in the cohomology
space. If the scalars $r_2$ and $r_3$ in the linear combination
are non-zero, then a meromorphic projective connection $R'$ and an
affine connections $T'$, both  holomorphic outside  $A$, can be
found such that $\ga=r_1 \ga_{S}^{(f)}+ r_2 \ga_{S,T'}^{(m)}+ r_3
\ga_{S,R'}^{(v)}$.
\end{theorem}

\subsection{Affine algebras}\label{SS:affine}
$ $

Let $\g$ be a reductive finite-dimensional Lie algebra. Above, we
introduced the current algebra $\gb$ together with its
almost-grading. In this subsection we study central extensions of
$\gb$. Given an invariant, symmetric bilinear form $\alpha(.,.)$,
i.e. a form obeying $\alpha([x,y],z)=\alpha(x,[y,z])$ we define as
generalization of the Kac-Moody algebras of affine type the {\it
higher genus (multi-point) affine Lie algebra}. We call it also a
{\it \KN\ algebra of affine type}. It is the Lie algebra based on
the vector space  $\gh=\C\oplus\gb$ equipped with  Lie structure
\begin{equation}\label{E:eaff}
[\widehat{x\otimes f},\widehat{y\otimes g}]= \widehat{[x,y]\otimes
(f g)}+\alpha(x,y)\cdot\gamma_{C_S}(f,g)\cdot t,\qquad
[\,t,\gh]=0\ ,
\end{equation}
where
\begin{equation}\label{E:gfun}
\gamma_{C_S}(f,g)=\cinc{C_S} fdg
\end{equation}
is the geometric cocycle for the function algebra obtained by
integration along a separating cycle $C_S$. We denote this central
extension by $\gh_{\alpha,S}$. It will depend on the bilinear form
$\alpha$. As usual we set $\widehat{x\otimes f}:=(0,x\otimes f)$.
The cocycle defining the central extension $\gh_{\a,S}$ is local.
Hence we can extend our almost-grading to the central extension by
setting $\deg t:=0$ and $\deg(\widehat{x\otimes A_{n,p}}):=n$.
Again we obtain a triangular decomposition
\begin{equation}\label{E:espaff}
\Gh_{\a,S}=\Gh_+\oplus\Gh_{(0)}\oplus\Gh_-\quad\text{with}\quad
\Gh_\pm\cong\G_\pm\quad\text{and}\quad\Gh_{(0)}=\G_{(0)}\oplus\C\cdot
t\ .
\end{equation}
The corresponding is true for the enlarged algebras. Among them,
$$
\gh^r:=\gh_-^{(1)}=\gb_-^{(1)}=\g\otimes \A_-^{(1)},\qquad
\gh_+^{*,ext}=\gb_+^*\oplus \C \;t=(\g\otimes \A_+^*)\oplus\C \;t\
.
$$
are of special interest.

Instead  integrating over a separating cycle in \refE{gfun} we
could take as integration path any other cycle $C$ and obtain in
this way other central extensions $\gh_{\a,C}$. In general
$\gh_{\alpha,C}$ will neither be equivalent nor isomorphic to
$\gh_{\alpha,S}$. In  addition, there is no a priori reason why
every cocycle defining a central extension of $\gb$ should be of
this type, i.e. should be obtained by choosing an invariant
symmetric bilinear form $\alpha$ and integrating the differential
$fdg$ over a cycle.

Before we can formulate the results needed in our context we have
to extend the definition of $\L$-invariance to $\gb$.
\begin{definition}\label{D:gblinv}
A cocycle $\ga$ of $\gb$ is called \emph{$\L$-invariant} if
\begin{equation}\label{E:glinv}
\ga(x(e\ldot g),y(h))+\ga(x(g),y(e\ldot h))=0,\qquad \forall
x,y\in\g,\ e\in\L,\ g,h\in\A.
\end{equation}
\end{definition}
The cocycles introduced above are obviously $\L$-invariant.
\begin{theorem}\label{T:simple}
\cite[Thm. 3.13, Cor. 3.14]{SchlMMJ} (a) Let $\g$ be a
finite-dimensional simple Lie algebra, then every local cocycle of
the current algebra $\gb=\g\otimes \A$ is cohomologous to a
cocycle given by
\begin{equation}\label{E:simple}
\gamma(x\otimes f,y\otimes g)=r\cdot
\frac{\b(x,y)}{2\pi\i}\int_{C_S}fdg,\quad \text{with}\ r\in\C,
\end{equation}
and with $\b$ the Cartan-Killing form of $\g$. In particular,
$\ga$ is cohomologous to a local and  $\L$-invariant cocycle.
\newline
(b) If the cocycle is already local and  $\L$-invariant, then it
coincides with the cocycle \refE{simple} with  $r\in\C$ suitable
chosen.
\newline
(c) For $\g$ simple, up to equivalence and rescaling of the
central element there is a unique non-trivial almost-graded
central extension $\gh$ of its higher genus multi-point current
algebra $\gb$. It is given by the cocycle \refE{simple}.
\end{theorem}
Next, let $\g$ be an arbitrary complex reductive
finite-dimensional Lie algebra, and
\begin{equation}\label{E:redecomp}
\g=\g_0\oplus\g_1\oplus\cdots\oplus\g_M
\end{equation}
be its decomposition into its abelian ideal $\g_0$ and simple
ideals $\g_1,\cdots,\g_M$. For the corresponding current algebra
we obtain
 $\gb=\gb_0\oplus\gb_1\oplus\cdots\oplus\gb_M$.
\begin{theorem}\label{T:red}
\cite[Thm. 3.20]{SchlMMJ} (a) Let $\g$ be a finite-dimensional
reductive Lie algebra. Given a cocycle $\ga$ for $\gb$ which is
local, and whose restriction to $\gb_0$ is $\L$-invariant, there
exists a symmetric invariant bilinear form $\a$ for $\g$ such that
$\ga$ is cohomologous to
\begin{equation}\label{E:ssimple}
 \gamma'_{\a,S}(x\otimes f,y\otimes g)=\frac{\a(x,y)}{2\pi\i}\int_{C_S}fdg.
\end{equation}
Vice versa, every such $\a$ determines a local cocycle.
\newline
(b) If the cocycle $\ga$ is already local and  $\L$-invariant for
the whole $\gb$, then it coincides with the cocycle $\ga'_{a,S}$.
\newline
(c)
\begin{equation}
\dim\H^2_{loc,\L}(\gb,\C)=\frac {n(n+1)}{2}+M.
\end{equation}
\end{theorem}
Here $\H^2_{loc,\L}(\gb,\C)$ denotes the subspace of those
cohomology classes which have a local and $\L$-invariant cocycle
as representative.


\subsection{Central extensions of $\D_\g$.}\label{SS:diffce}
$ $

Recall that  $\D_\g=\gb\oplus\L$ with $[e,xA]=x(e.A)$, see
\refE{comm}. There is the short exact sequence of Lie algebras
\begin{equation}\label{E:dsequ}
\begin{CD}
0@>>>\gb@>i_1>>\Do_\g@>p_2>>\L@>>>0.
\end{CD}
\end{equation}
First note again that by using the almost-grading of $\L$ and
$\gb$ and by using the fact that $\A$ is an almost-graded
$\L$-module the algebra $\Do_\g$ is an  almost-graded algebra.
Again, in \cite{SchlMMJ} local cocycles $\ga$ and central
extensions of this algebra are considered. By restricting a local
cocycle for $\D_\g$ to the subalgebra $\gb$ one obtains a local
cocycle for $\gb$. In the semi-simple case we obtain
\begin{theorem}\label{T:dsimple}
(a) Let $\g$ be a semi-simple Lie algebra and $\ga$ a local
cocycle of $\D_\g$. Then there exists a symmetric invariant
bilinear form $\a$ for $\g$ such that $\ga$ is cohomologous to a
linear combination of the local cocycle $\ga_{\a,S}$ given by
\refE{ssimple}  and of the local cocycle $\ga_{S,R}^{(v)}$
\refE{ecva} for $C=C_S$ of the vector field algebra $\L$.
\newline
(b) If $\g$ is a simple Lie algebra, then $\ga_{\a,S}$ is a
multiple of the standard cocycle \refE{simple} for $\g$.
\newline
(c) $\dim\H^2_{loc}(\D_\g,\C)=M+1$, where $M$ is the number of
simple summands of $\g$.
\end{theorem}
In the reductive case it turns out that a cocycle of $\D_\g$
restricted to the abelian summand $\gb_0$ is $\L$-invariant. In
generalization of the mixing cocycle for $\D$ we obtain for every
linear form $\phi\in\g^*$ which vanishes on
$\g':=[\g,\g]=\g_1\oplus\cdots \oplus\g_M$ a local cocycle given
by
\begin{equation}\label{E:amix}
\ga_{\phi,S}(e,x(g)):=\frac {\phi(x)}{2\pi\i} \int_{C_S}
\left(\tilde e\cdot g''+T\cdot (\tilde e\cdot
 g')\right)dz,
\end{equation}
Here again $T$ is a meromorphic connection with poles only at the
points in $A$.
\begin{theorem}\label{T:dgred}
\cite[Thm. 4.11]{SchlMMJ} (a) Let $\g$ be a finite-dimensional
reductive Lie algebra. For every local cocycle $\ga$ for $\D_\g$
there exists a symmetric invariant bilinear form $\a$ of $\g$, and
a linear form $\phi$ of $\g$ which vanishes on $\g'$, such that
$\ga$ is cohomologous to
\begin{equation}\label{E:gred}
\ga'=\ga_{\a,S}+\ga_{\phi,S}+r\ga_{S,R}^{(v)},
\end{equation}
with $r\in\C$ and a current algebra cocycle $\ga_{\a,S}$ given by
\refE{ssimple},  a mixing cocycle $\ga_{\phi,S}$ given by
\refE{amix} and the vector field cocycle $\ga_{S,R}^{(v)}$ given
by \refE{ecva}. Vice versa, any such $\a$, $\phi$, $r\in\C$
determine a local cocycle.
\newline
(b) The space of local cocycles $\H_{loc}^2(\D_\g,\C)$ is $\dfrac
{n(n+1)}{2}+n+M+1$ dimensional.
\end{theorem}

\subsection{Local cocycles for $\sln(n)$ and $\gl(n)$}
$ $

Consider $\sln(n)$, the Lie algebra of trace-less complex $n\times
n$ matrices. Up to multiplication with a scalar the Cartan-Killing
form $\b(x,y)=\tr(xy)$ is the unique symmetric invariant bilinear
form. From the Theorems \ref{T:simple} and \ref{T:dsimple} follows
\begin{proposition}\label{P:slcoc}
(a) Every local cocycle for the current algebra
$\overline{\sln}(n)$ is cohomologous to
\begin{equation}\label{E:slcoc}
\ga(x(g),y(h))=r\cdot \frac {\tr(xy)}{2\pi\i}\int_{C_S} gdh,\quad
r\in\C.
\end{equation}
(b) Every $\L$-invariant local cocycle equals the cocycle
\refE{slcoc} with a suitable $r$.
\newline
(c) Every local cocycle for the differential operator algebra
$\Do_{\sln(n)}$ is cohomologous to a linear combination of
\refE{slcoc} and the standard local cocycle $\ga_{S,R}^{(v)}$ for
the vector field algebra. In particular, there exist no cocycles
of pure mixing type.
\end{proposition}
Next, we deal with $\gl(n)$, the Lie algebra of all complex
$n\times n$-matrices. Recall that $\gl(n)$ can be written as Lie
algebra direct sum $\gl(n)=\mathfrak{s}(n)\oplus \sln(n)\cong\C
\oplus\sln(n)$. Here $\mathfrak{s}(n)$ denotes the $n\times n$
scalar matrices. This decomposition is the decomposition as
reductive Lie algebra into its abelian and semi-simple summand.

The space of symmetric invariant bilinear forms for  $\gl(n)$ is
two-dimensional. A basis is given by the forms
\begin{equation}
\a_1(x,y)=\tr(xy),\quad\text{and}\quad \a_2(x,y)=\tr(x)\tr(y).
\end{equation}
The form $\a_1$ is the ``standard'' extension for the
Cartan-Killing form for $\sln(n)$ to $\gl(n)$ and is also $\gl(n)$
invariant. From the Theorems \ref{T:red} and \ref{T:dgred} follows
\begin{proposition}\label{P:glcoc}
(a) A cocycle $\gamma$ for $\overline{\gl}(n)$ is local and
restricted to $\overline{\mathfrak{s}(n)}$ is $\L$-invariant if
and only if it is cohomologous to a linear combination of the
following two cocycles
\begin{equation}\label{E:glns}
\ga_1(x(g),y(h))=\frac{\tr(xy)}{2\pi\i} \int_{C_S} gdh,\qquad
\ga_2(x(g),y(h))=\frac{\tr(x)\tr(y)}{2\pi\i} \int_{C_S} gdh.
\end{equation}
(b) If the local cocycle $\ga$  restricted to $\overline{\gl}(n)$
is $\L$-invariant then $\ga$ is already equal to this linear
combination of the cocycles \refE{glns}.
\end{proposition}
\begin{proposition}\label{P:dgln}
(a) Every local cocycle $\ga$ for  $\Do_{\gl(n)}$ is cohomologous
to a linear combination of the cocycles $\ga_1$ and $\ga_2$ of
\refE{glns}, of the mixing  cocycle
\begin{equation}
\ga_{3,T}(e,x(g))=\frac {\tr(x)}{2\pi\i}\int_{C_S} \big(\tilde e
g''+T\tilde e g'\big)dz,
\end{equation}
and of the standard local cocycle $\ga_{S,R}^{(v)}$ for the vector
field algebra, i.e.
\begin{equation}\label{E:decompc}
\gamma=r_1\gamma_1+r_2\gamma_2
+r_3\ga_{3,T}+r_4\ga_{S,R}^{(v)}+\text{coboundary},
\end{equation}
with suitable $r_1,r_2,r_3,r_4\in\C$.
\newline
(b) If the cocycle $\ga$ is local and restricted to
$\overline{\gl}(n)$ is $\L$-invariant, and $r_3,r_4\ne 0$ then
there exist an affine connection $T$ and a projective connection
$R$ holomorphic outside $A$ such that in the combination
\refE{decompc} there will be no additional coboundary.
\newline
(c) $\dim\H^2_{loc}(\Do_{\gl(n)},\C)=4$.
\end{proposition}
It turns out that for $\gamma$ as a cocycle of the differential
operator algebra its restriction to
 $\overline{\gl}(n)$ is
cohomologous to an $\L$-invariant cocycle. Moreover
$\ga_{|\mathfrak{s}(n)}$ is already $\L$-invariant, see
\cite[Prop. 4.10]{SchlMMJ}.

The mixing type cocycle for ${\gl}(n)$ was also studied in
\cite{ShMMJ} .
\section{Representations  of the multi-point \KN\ algebras}
\label{S:kzrep}
%
\subsection{Projective $\D_\g$-modules}${ }$\label{SS:pro}

Let $\V$ be an arbitrary Lie algebra, $V$ a vector space and $\pi:
\V \to End(V)$ a linear map. The space $V$ is called a {\it
projective $\V$-module} if for all pairs $f,g\in\V$ there exists
$\ga(f,g)\in\C$ such that
\begin{equation}\label{E:pro}
 \pi([f,g])=[\pi(f),\pi(g)]+\ga(f,g)\cdot id.
\end{equation}
In this case $\pi$ is also called a projective action. Often we
write simply $f\,v$ instead  $\pi(f)(v)$.

If $V$ is a projective $\V$-module then $\ga$ is necessarily a Lie
algebra 2-cocycle. Via Equation \refE{centext} the cocycle,  hence
the projective action $\pi$, defines a central extension
${\widehat\V}_\ga={\widehat\V}_\pi$ of $\V$. Obviously, by
defining $\pi(t)=id$, the projective action can be extended to a
honest Lie action of  ${\widehat\V}_\ga$ on $V$. Usually, we
suppress the index $\pi$ if it is clear from the context. One
should  keep in mind that the cocycle defining $\widehat\V$ and
(without further restriction on the projective action) also the
equivalence and even isomorphy class of the central extension will
depend on the projective module action  given. In certain cases
considered below we will be able to identify  the cocycle.

Let $\V$ be an almost-graded Lie algebra. A projective $\V$-module
$V$ is called {\it admissible} if for each $v\in V$ there exists a
$k(v)\in\N$ such that for all  $f\in\V$ with $\deg f\ge k(v)$ we
have $fv=0$. By a projective almost-graded module $V$ over an
almost-graded algebra $\V$ we understand an almost-graded module
structure on $V$ with respect to the projective action of $\V$
such that the corresponding cocycle is local (see \refD{local}).
We call an almost-graded module $V$ a {\it highest weight (vacuum)
module} if an element $\vac\in V$ exists such that $f\vac =0$ for
every $f$ with   $\deg f>0$ and $V=U(\V)\vac$ where $U(\cdot)$
denotes the universal enveloping algebra. The element $\vac$ is
called a vacuum vector. Obviously, vacuum modules are admissible.

Our Lie algebra will be mainly $\D_\g$.
 We will assume
$V$ to be an admissible projective almost-graded $\D_\g$-module
(we will write {\it projective representation of $\D_\g$} as
well). Recall that via the  embedding \refE{dsequ} the space $V$
is also a projective $\gb$-module. The most important examples of
such modules, which, moreover, are generated by vacuum vectors,
are  the fermion representations which will be discussed in the
next subsection.

In fact, we will be mainly interested in the case
$\g=\frak{gl}(n)$. In particular, the fermion modules (of fixed
charge) are irreducible in this case. Observe that for the
constructions of the following \refS{kzgen} we need only a
(projective) $\gb$-module structure on $V$, and we will assume the
irreducibility of this module.

\subsection{Fermion representations}\label{SS:ferm}${ }$

In this paragraph, we briefly outline the construction of a
projective fermion $\D_\g$-module \cite{rShf,ShMMJ,ShN65}. It uses
the Krichever-Novikov bases in the spaces of sections of
holomorphic vector bundles as an important ingredient. These bases
are introduced in \cite{rKNR2p} for the two-point case, and,
combining the approaches of \cite{rKNR2p} and \refSS{knb} (going
back to \cite{rSLa,rSLb}) generalized in \cite{ShN65} to the
multi-point case. We refer to the cited works for the details.
Here we only mention that these bases are given by asymptotic
behavior of their elements at the marked points and their
properties are similar to those described in \refSS{knb}. In
particular, these bases are almost graded with respect to the
action of Krichever-Novikov algebras.

Consider a holomorphic bundle $F$ on $\Sigma$ of rank $r$ and
degree $g\cdot r$ ($g$ is the genus of our Riemann surface
$\Sigma$). Let $\Ga(F)$ denote the space of meromorphic sections
of $F$ holomorphic except at $P_1,\ldots,P_N,\Pif$. Let $\tau$ be
a finite-dimensional representation of $\g$ with representation
space $V_\tau$. Set $\Ga_{F,\tau}:=\Ga(F)\otimes V_\tau$.

We define a $\gb$-action on $\Ga_{F,\tau}$ as follows:
\begin{equation}\label{E:gbact}
       (x\otimes A)(s\otimes v)=(A\cdot s)\otimes \tau(x)v\quad
           \text{for all}\quad
        x\in\g, A\in\A, s\in\Ga(F),v\in V_\tau .
\end{equation}

In order to define an $\L$-action on $\Ga(F)$ we choose a
meromorphic (therefore flat) connection $\nabla$ on $F$ which has
logarithmic singularities at $P_1,\ldots,P_N$ and $\Pif$ (see
\cite{ShMMJ} for more details). By flatness,
$\nabla_{[e,f]}=[\nabla_e,\nabla_f]$ for all $e,f\in\L$. Hence,
$\nabla$ defines a representation of $\L$ in $\Ga(F)$.

{}From  the definition of a connection, for any $s\in\Ga(F)$,
$e\in\L$ and $A\in\A$ we have $\nabla_e(As) = (e\ldot
A)s+A\nabla_es$ where $e\ldot A$ is the Lie derivative. Hence,
$[\nabla_e,A]=e\ldot A$, i.e.  the mapping
$e+A\rightarrow\nabla_e+A$ gives rise to a representation of $\D$
in $\Ga(F)$.

Define the corresponding $\L$-action on $\Ga_{F,\tau}$ by
\begin{equation}\label{E:lac}
       e(s\otimes v)=\nabla_es\otimes v \quad\text{for all}\quad
        e\in\L , s\in\Ga(F),v\in V_\tau .
\end{equation}
It can be verified directly that \refE{gbact} and \refE{lac} give
a representation of $\D_\g$ in $\Ga_{F,\tau}$.

\medskip
Choose a \KN\ basis in $\Ga(F)$ \cite{ShN65} and a weight basis
$\{v_i|1\le i\le \dim V_\tau\}$ in $V_\tau$. Each basis element in
$\Ga(F)$ is determined by its {\it degree} $n\in{\mathbb Z}$, the
number $p$ of a marked point: $1\le p\le N$,  and  an integer $j$:
$0\le j\le r-1$. Denote the element corresponding to a triple
$(n,p,j)$ by $\psi_{n,p,j}$.  Introduce
$\psi_{n,p,j}^i=\psi_{n,p,j}\otimes v_i$. Enumerate the elements
$\psi_{n,p,j}^i$ linearly in ascending lexicographical  order of
the quadruples $n,p,j,i$. In this way we set
$\psi_M=\psi_{n,p,j}^i$ where $M=M(n,p,j,i)\in\Z$ is in linear
order.

\begin{lemma}\label{L:agda}\cite{ShN65}
With respect to the index M, the module $\Ga_{F,\tau}$ is an
almost-graded $\D_\g$-module.
\end{lemma}
The proof of the almost-gradedness is similar to that in the
two-point case \cite{rShf,ShMMJ}.

\medskip
The final step of the construction of the fermion representation
corresponding to the pair $(F,\tau)$ is passing to the space of
the semi-infinite monomials on $\Ga_{F,\tau}$.

Consider the vector space $\Hwft$ generated over $\mathbb C$ by
the formal expressions ({\it semi-infinite monomials}) of the form
$\Phi=\psi_{N_0}\wedge\psi_{N_1} \wedge\ldots$, where the
$\psi_{N_i}$ are the above introduced basis elements of
$\Ga_{F,\tau}$, the indices are strictly increasing, i.e.
$N_0<N_1<\ldots$, and for all $k$ sufficiently large    $N_k=k+m$
for a suitable $m$ (depending on the monomial).
 Following \cite{rKaRa}
we call $m$ {\it the charge} of the monomial. For a monomial
$\Phi$ of charge $m$ the degree of $\Phi$ is defined as follows:
\begin{equation}\label{E:mdeg}
    \deg\Phi =\sum\limits_{k=0}^\infty (N_k-k-m).
\end{equation}
Observe that there is an arbitrariness in the enumeration of the
$\psi_{n,p,j}^i$'s for a fixed $n$; the just defined degree of a
monomial does not depend on this arbitrariness.

We want to extend the action of $\D_\g$ on $\Ga_{F,\tau}$ to
$\Hwft$. Here, we briefly outline the construction; see
\cite{rShf,ShMMJ} for details. Assuming the $\D_\g$-action on
$\psi_{N_i}$'s to be known, apply a basis element of $\D_\g$ to
$\Phi$ by the Leibnitz rule. If, in the process, a monomial
containing the same $\psi_N$ in different positions occurs, it is
set to zero. If a pair $\psi_N\wedge\psi_{N'}$ in a wrong order
($N>N'$) occurs then it should be transposed and the sign before
the corresponding monomial changes. This is  done until all
entries are in the strictly increasing order. Due to the
almost-gradedness of $\D_\g$-action on $\Ga_{F,\tau}$ and the
above mentioned stabilization ($N_k=k+m$, $k\sim\infty$), the
result of the above steps is well defined for all basis elements
 except for the finite number. For
those, apply the standard \emph{process of regularization}
\cite{rKaRa,rKNFb,rKNFc,rSDiss}. In this way the
 action on $\Ga_{F,\tau}$ can be extended to  $\Hwft$ as a
projective Lie algebra action.

Let $\Hwftm$ be the subspace of $\Hwft$ generated by the
semi-infinite monomials of charge $m$. These subspaces are
invariant under the projective action of $\D_\g$. This follows, as
in the classical situation, from the fact that after the action
of $\D_\g$ the resulting monomials will have the same ``length''
and the same ``tail'' as the monomial one has  started with.
Hence, for every $m$ the space $\Hwftm$ is itself a projective
$\D_\g$-module, and $\Hwft=\bigoplus_{m\in\Z}\Hwftm$ as projective
$\D_\g$-module.
 We call the modules $\Hwft$,$\Hwftm$
(projective) {\it fermion representations}.

\begin{proposition}\label{P:fermalm}
Let $\Hwftm$ be the submodule of $\Hwft$ of charge $m$.
\newline
(a) With respect to the degree \refE{mdeg} the homogeneous
subspaces $(\Hwftm)_k$ of degree $k$ are finite-dimensional. If
$k>0$ then $(\Hwftm)_k=0$.
\newline
(b) The cocycle $\gamma$ for  $\D_\g$  defined by the projective
representation is local. It is bounded from above by zero.
\newline
(c) The module $\Hwftm$ is an almost-graded projective
$\D_\g$-module.
\end{proposition}
\begin{proof}
(a) From the very definition of the degree it follows that the
degree $k$ is always $\le 0$. Notice that all summands in
\refE{mdeg} are $\le 0$. Hence for a given $k$ only finitely many
semi-infinite monomials of charge $m$ can realize this $k$.
\newline
(b) and (c) follow as in the two-point case, \cite{rShf,ShMMJ}.
\end{proof}
From this proposition and the classification results of
\refT{dgred} we obtain
\begin{proposition}
The module   $\Hwftm$ is a Lie module over a certain central
extension $\widehat{\D}_\g$ defined via a local cocycle of the
form \refE{gred}.
\end{proposition}
Recall that by restricting the action to $\gb$ we obtain an
admissible projective representation and hence an almost graded
central extension $\gh$ of $\gb$. In particular, for $\g$
reductive this restriction is cohomologous to an $\L$-invariant
cocycle, i.e. it is cohomologous to a geometric cocycle of the
type \refE{eaff} with a suitable invariant symmetric bilinear form
$\alpha$ (see \refT{red}). For the special case $\g=\gl(n)$ the
cocycle classes are given in \refP{glcoc} and \ref{P:dgln}.

Under an admissible representation $V$ of $\gh$ we understand a
representation admissible with respect to the almost-grading (in
the sense introduced above) as projective representation of $\gb$
in which the central element $t$ operates as a scalar $\ce\cdot
id$, $\ce\in{\mathbb C}$. The number $\ce$ is called the {\it
level} of the $\gh$-module $V$. It follows immediately
\begin{proposition}
 $\Hwftm$ carries an admissible representation of $\gh$.
\end{proposition}

\subsection{Sugawara representation}${ }$

First let $\g$  be any finite dimensional reductive Lie algebra.
We fix an invariant symmetric bilinear form $\alpha$ on $\g$.
Starting from this section we have to assume that $\alpha$ is
non-degenerate. By $\gh$ we denote the standard central extension
(depending on the bilinear form $\alpha$) as introduced in
\refSS{affine} (see Equations \refE{eaff} and \refE{gfun})
together with its almost-grading. Let $V$ be any admissible
representation of level $c$.

If $\g$ is abelian or simple then each admissible representation
of $\gh$ of non-critical level, i.e. for a level which is not the
negative of the dual Coxeter number in the simple case, or a level
$\ne 0$ in the abelian case,  the (affine) Sugawara construction
yields a projective representation of the \KN\ vector field
algebra $\L$. This representation is called the {\it Sugawara
representation}.  For the two-point case, the abelian version of
this construction was introduced in \cite{rKNFb}. The nonabelian
case was later considered in \cite{rBono}, \cite{rSSS}.  In
\cite{rSSS} also the multi-point version was given.
 Observe that
every positive level is non-critical. For an arbitrary complex
reductive Lie algebra $\g$ the Sugawara representation is defined
as a certain linear combination of Sugawara representations of its
simple ideals (V.Kac \cite{KacB}, \cite[Lecture 10]{rKaRa}).

Due to its importance in our context, we have to describe the
construction in more detail. For any $u\in\g$, $A\in \A$ we denote
by $\ u(A)\ $ the operator in $V$ corresponding to $u\otimes A$.
We also denote an element of the form $u(A_{n,p})$  by $\ u(n,p)\
$. We choose a basis $\ u_i,\ i=1,\ldots,\dim\g\ $ of $\g$ and the
corresponding dual basis $\ u^i,\  i=1,\ldots,\dim\g$ with respect
to the form $\a$.  The {\it Casimir element} $\
\Omega^0=\sum_{i=1}^{\dim \g} u_iu^i\ $ of the universal
enveloping algebra  $U(\g)$ is independent of the choice of the
basis. To simplify notation, we denote $\sum_i u_i(n,p)u^i(m,q)$
by $u(n,p)u(m,q)$.

We define the higher genus {\it Sugawara operator} (also called
{\it Segal operator} or {\it energy-momentum tensor}) as
\begin{equation}\label{E:suga}
T(P):=\frac 12\sum_{n,m}\sum_{p,s}
\nord{u(n,p)u(m,s)}\w^{n,p}(P)\w^{m,s}(P)\ .
\end{equation}
By $\ \nord{....}\ $ we denote some normal ordering. In this
section the summation indices $n,m$ run over $\Z$, and $p,s$ over
$\{1,\ldots,N\}$. The precise form of the normal ordering is of no
importance here.  As an example we may take the following
``standard normal ordering'' ($x,y\in\g$)
\begin{equation}\label{E:normst}
\nord{x(n,p)y(m,r)}\ :=
\begin{cases} x(n,p)y(m,r),&n\le m
                                           \\
                            y(m,r)x(n,p),&n>m\ .
                      \end{cases}
\end{equation}

The expression $T(P)$ can be considered as a formal series of
quadratic differentials in the variable $P$ with operator-valued
coefficients. Expanding it over the basis $\Omega^{k,r}$ of the
quadratic differentials we obtain
\begin{equation}\label{E:sugb}
  T(P)=\sum_k\sum_r L_{k,r}\cdot\Omega^{k,r}(P)\ ,
\end{equation}
with
\begin{equation}\label{E:sugc}
\begin{gathered}
  L_{k,r}=\cins T(P)e_{k,r}(P)=\frac 12\sum_{n,m}\sum_{p,s}
  \nord{u(n,p)u(m,s)}l_{(k,r)}^{(n,p)(m,s)},\\
  \text{where}\qquad
  \ l_{(k,r)}^{(n,p)(m,s)}:=\cins \w^{n,p}(P)\w^{m,s}(P)e_{k,r}(P)\ .
\end{gathered}
\end{equation}
Formally, the operators $L_{k,r}$ are  infinite double sums. But
for given $k$ and $m$, the coefficient $l_{(k,r)}^{(n,p)(m,s)}$
will be non-zero only for finitely many $n$. This can be seen by
checking the residues of the elements appearing under the
integral. After applying the remaining infinite sum to a fixed
element $v\in V$, by the normal ordering and admissibility of the
representation only finitely many of the operators will operate
non-trivially on this element.

The following theorem is proved in \cite{rSSS}.
\begin{theorem}\label{T:suga}
Let $\g$ be a finite dimensional either abelian or simple Lie
algebra and $2\k$ be the eigenvalue of its Casimir operator in the
adjoint representation. Let $\alpha$ be the normalized Cartan
Killing form in the simple case or any non-degenerate bilinear
form in the abelian case, and $\gh$ be the corresponding central
extension. Let $V$ be  an admissible almost-graded $\gh$-module of
level $\ce$. If $\ce+\k\ne 0$ then the rescaled ``modes''
\begin{equation}\label{E:sugm}
L_{k,r}^*=\frac {-1}{2(\ce+\k)}\sum_{n,m}\sum_{p,s}
 \nord{u(n,p)u(m,s)}l_{(k,r)}^{(n,p)(m,s)}
\ ,
\end{equation}
of the Sugawara operator are well-defined operators on $V$ and
define an admissible projective representation of $\L$. The
corresponding cocycle for the vector field algebra $\L$ is a local
cocycle.
\end{theorem}
\begin{remark}
By the locality of the cocycle and in view of  \refT{unique} the
cocycle is  a scalar multiple of \refE{ecva} with $C=C_S$. In
particular, the central extension $\Lh$ for which the Sugawara
representation is a honest  representation, is fixed up to
isomorphy.
\end{remark}
\begin{proposition}\label{P:sugalm}
The module $V$ is an almost-graded $\Lh$-module under the Sugawara
action.
\end{proposition}
\begin{proof}
We have to show that there exist constants $M_1,M_2$ such that for
every given  homogeneous element $\psi_s\in V$ of degree $s$, and
every $k$ we have
\begin{equation}\label{E:krt}
k+s+M_1\le \deg(L_{k,r}^*\psi_s)\le k+s+M_2.
\end{equation}
Starting from the description \refE{sugm} we consider first the
coefficients $l_{(k,r)}^{(n,p)(m,s)}$, which are given as
integrals \refE{sugc}. The integral could only be non-vanishing if
the integrand has poles at the points in $I$ and at the point
$\Pif$. Using the explicit formulas \refE{ordfn} and \refE{ordi}
for the orders at these points we obtain
\begin{equation}\label{E:kgN}
k\le n+m\le k+C(g,N),
\end{equation}
with a rational constant $C(g,N)\ge 0$ only depending on the genus
$g$ and the number of points $N$. By the almost-gradedness of $V$
as  $\gh$-module, there exist constants $c_1$ and $c_2$ such that
for all $m,r,s$
\begin{equation}
m+s+c_1\le \deg(u(m,r)\psi_s)\le m+s+c_2,
\end{equation}
if the element $u(m,r)\psi_s\ne 0$. Hence,
\begin{equation}
n+m+s+2c_1\le \deg(\nord{u(n,p)u(m,r)}\psi_s)\le n+m+s+2c_2,
\end{equation}
if the element $\nord{u(n,p)u(m,r)}\psi_s\ne 0$. Using \refE{kgN}
we obtain Equation \refE{krt} if we set $M_1=2c_1$ and
$M_2=2c_2+C(g,N)$.
\end{proof}
We call the $L_{k,r}^*$, resp. the $L_{k,r}$ the Sugawara
operators too. For $e=\sum_{n,p}a_{n,p}e_{n,p}\in\L$
($a_{n,p}\in{\mathbb C}$) we set $T[e]=\sum_{n,p}a_{n,p}L_{n,p}^*$
and obtain the projective representation $T$ of $\L$. It is called
the projective {\it Sugawara representation} of the Lie algebra
$\L$ corresponding to the given admissible representation $V$ of
$\gh$.

By the Krichever-Novikov duality the Sugawara operator $T[e]$
assigned to the vector field
 $e\in\L$ can equivalently be given as
\begin{equation}\label{E:suge}
T[e]=\frac {-1}{\ce+\ka}\cdot \cins T(P)e(P).
\end{equation}

Let  $\g$ be a reductive Lie algebra with decomposition
\refE{redecomp}. The elements $x\in\g$ can be decomposed as
$x=\sum_{i=0}^n x_i$, with $x_i\in\g_i$. Let $\alpha$ be a
symmetric invariant  bilinear form for $\g$. With respect to this
decomposition from the invariance follows $\alpha(x_i,x_j)=0$ for
$i\ne j$, i.e. the decomposition is an orthogonal decomposition.
By restricting $\alpha$  to $\g_i$ we obtain a symmetric invariant
bilinear form on $\g_i$, hence a multiple of the Cartan-Killing
form (resp. an arbitrary symmetric bilinear form on $\g_0$). A
non-degenerate form is called {\it normalized} if the restrictions
of $\alpha$ for the simple summands is equal to the Cartan-Killing
form.

Let $\ga$ be a local cocycle for $\gb$ for $\g$ reductive,
$A,B\in\A$, and $x,y\in\g$. It follows from the cocycle condition
(see \cite[Lemma 3.11]{SchlMMJ}) that  $\ga(x_i A,y_j B)=0$ for
$i\ne j$ (with the same decomposition for $y$ as for $x$ above).
This implies
\begin{equation}
\ga(x A,y B)=\sum_{i=0}^M\ga(x_i A,y_i B).
\end{equation}
Given an admissible  representation $V$ of $\gh$  we obtain
representations of $\gh_i$ on the same space $V$. We can define
the individual (rescaled) Sugawara operators $T_k[e]$,
$k=0,1,\ldots ,M$. We set
\begin{equation}
T[e]:=\sum_{k=0}^M T_k[e].
\end{equation}
The following lemma expresses a fundamental property of the
Sugawara representation. It was shown in \cite{rSSS} for the
abelian and simple case. Here we will show how to extend the
result to the general reductive case.
 \begin{lemma}\label{L:fund}
 Let $\g$ be a reductive Lie algebra with a chosen normalized form
$\alpha$, and $T[e]$ for every $e\in\L$ the operator  as defined
above. Then $T$ defines a representation of the  centrally
extended vector field algebra $\Lh$ and for any $x\in\g$,
$A\in\A$, $e\in\L$, we have
\begin{equation}\label{E:frel}
    [T[e],x(A)]=x(e\ldot A).
\end{equation}
\end{lemma}
\begin{proof}
Given $x,y\in\g$ denote by $x_i$ and $y_j$ its components as
above. Then $[x_i(A),y_j(B)]=[x_i,y_j](AB)+\ga(x_i A,y_j
B)\,c\cdot id$. In particular, if $i\ne j$ we obtain that $x_i$
and $x_j$ commute and that the cocycle vanishes. Hence these
elements commute also. Let $T_k$,  be the Sugawara representation
corresponding to the representation $V$ of $\widehat{\g_k}$. Since
the operators of the representations $T_i$ are expressed via
$x_i(A)$'s and operators of $T_j$ via $x_j(A)$'s, $T_i$ and $T_j$
commute for $i\ne j$. Hence, $T$ is a representation of $\Lh$.
Moreover, for $x_i\in\g_i$, $x_j\in\g_j$, $A\in\A$, $e\in\L$ we
have  $[T_j[e],x_i(A)]=0$. For $\g$ simple or abelian, \refE{frel}
is shown in \cite{rSSS} (see also \cite{rSSpt}). Hence,
\[ [T_k[e],x_k(A)]=x_k(e\ldot A),\quad k=0,1,\ldots,M.
\]
This implies
\[
    [T[e],x(A)]
     =[\sum_{k=0}^M T_k[e],\sum_{i=0}^M  x_i(A]
     = \sum_{i=0}^M x_i(e.A)= x(e.A).
\]
 \end{proof}

\medskip
Having only one simple summand for $\g=\gl(n)$ the condition that
the form $\alpha$ is normalized on its simple summands can always
be achieved by rescaling the level. In case of the fermion
representations for $\gl(n)$ we obtain that for these
representations relation \refE{frel}  is true.

\section{Moduli of curves with marked points,
 conformal blocks and projective flat connection}\label{S:kzgen}

\subsection{Moduli space $\MgNekp$ and the sheaf of
conformal blocks}${ }$

In \cite {rSSpt} we described the moduli spaces of curves which
typically  occur in 2d conformal field theories. Here we will
slightly extend the definitions introduced there. We denote by
$\MgNekp$ the moduli space of smooth projective curves of genus
$g$ (over $\C$) with $N+1$ ordered distinct marked points and
fixed $k$-jets of local coordinates at the first $N$ points and a
fixed $p$-jet of a local coordinate at the last point. The
elements of $\MgNekp$ are given as
\begin{equation}\label{E:point}
   \bt^{(k,p)}=[\Sigma,P_1,\ldots,P_N,P_\infty,
   z_1^{(k)}\ldots,z_N^{(k)},z_\infty^{(p)}]\ ,
 \end{equation}
 where $\Sigma$ is a smooth projective  curve of genus $g$,
$P_i$, $\iNi$ are  distinct points on $\Sigma$,
 $z_i$ is a coordinate at $P_i$
with $z_i(P_i)=0$, and $z_i^{(l)}$ is a $l$-jet of $z_i$ ($l\in
\N_0$). Here $[..]$ denotes an equivalence class of such tuples in
the following sense. Two tuples representing $\bt^{(k,p)}$ and
${\bt}^{(k,p)}{}'$ are equivalent if there exists an algebraic
isomorphism $\phi:\Sigma\to \Sigma'$ with $\phi(P_i)=P_i'$ for
$i=1,\ldots,N,\infty$ such that after the identification via
$\phi$ we have
\begin{equation}
  z_i'=z_i+O(z_i^{k+1}), \quad \iN\ \quad\text{and}\quad
  z_\infty'=z_\infty+O(z_i^{p+1}).
\end{equation}
For the following two special cases we introduce the same notation
as in \cite{rSSpt}: $\MgNe= \MgNezz$, and $\MgNep=\MgNeee$ . By
forgetting either coordinates or higher order jets we obtain
natural projections
\begin{equation}
\MgNeep\to\MgNe, \qquad \MgNekp\to \MgNekkpp
\end{equation}
for any $k'\le k$ and $p'\le p$. In this article (as well as in
the previous one \cite {rSSpt})
 we are only dealing
with the local situation in the neighborhood of a moduli point
corresponding to a generic curve $\Sigma$ with a generic marking
$(P_1,P_2,\ldots, P_N,\Pif)$. Let $\Wt\subseteq \MgNe$ be an open
subset  around such a generic  point $\tilde b=[\mpt]$. A generic
curve of $g\ge 2$ admits no nontrivial infinitesimal automorphism,
and we may assume that there exists a universal family of curves
with marked points over $\Wt$. In particular, this says that
there is a proper, flat family of smooth curves over $\Wt$
\begin{equation}\label{E:unim}
  \pi :\Cal U\to  \Wt\ ,
\end{equation}
such that for the points $\tilde b=[\mpt]\in\Wt$ we have $\
\pi^{-1}(\tilde b)= \Sigma$ and that the sections defined as
\begin{equation}\label{E:sect}
\sigma_i:\Wt\to\Cal U,\quad\sigma_i(\tilde b)=P_i,\quad \iN,\infty
\end{equation}
are holomorphic. For more background information, see \cite[Sect.
1.2, Sect. 1.3]{rUcft}, in particular Thm. 1.2.9 of \cite{rUcft}.

If we ``forget'' the last point $\Pinf$ we obtain maps
\begin{equation}
\MgNe\to\MgN,\quad \MgNezp\to\MgN,\quad \MgNekp\to\MgNk.
\end{equation}
Let us fix a holomorphic section $\sinf$ of the universal family
of curves (without marking). In particular, for every curve there
is a point chosen in a manner depending analytically on the
moduli.  (Recall, we are only dealing with the local and generic
situation.) The analytic subset
\begin{equation}\label{E:msidw}
W':=\{\tilde b=[\mpt]\mid
\Pif=\sinf([\Sigma])\}\quad\subseteq\quad \Wt
\end{equation}
can be identified with an open subset $W$ of $\MgN$ via
\begin{equation}\label{E:msid}
\tilde b=[ (\Sigma,P_1,P_2,\ldots, P_N,\sinf([\Sigma]))] \ \to\
b=[ (\Sigma,P_1,P_2,\ldots, P_N)]\ .
\end{equation}
By genericity, the map is one-to-one.

By choosing not only a section $\sinf$ but also a $p$-th order
infinitesimal neighborhood of this section we even get an
identification of the open subset $W$ of $\MgN$ with an analytic
subset $W^{\prime,(p)}$ of $\Wt^{(0,p)}$ of $\MgNezp$.  It is
defined in a similar way as $W'$.

All these considerations can be extended to the case where we
allow infinite jets of local coordinates at $\Pinf$. We obtain
then the moduli space $\MgNeki$.

\medskip

At generic points, the moduli spaces $\MgNekp$ are smooth. Denote
by $S$ the divisor $S=\sum_{i=1}^NP_i$ on $\Sigma$. The tangent
space $\T_{\bt^{(1,p)}}\MgNeep$ can be identified with the
cohomology space $\H^1(M,T_M(-2S-(p+1)\Pinf))$. As in Prop. 4.4
and Thm. 4.5 of \cite{rSSpt} we obtain that there exists a
surjective linear map from the Krichever-Novikov vector field
algebra $\L$ to the cohomology space
\begin{equation}
\theta=\theta_{p}:\L\to \H^1(M,T_M(-2S-(p+1)\Pinf)
\end{equation}
such that $\theta$ restricted to the following subspaces
 gives  isomorphisms
\begin{equation}\label{E:modiso}
\begin{aligned}
\L_0\oplus\L_{-1}\oplus\L_{(0)}^{(p)}&\cong
\H^1(M,T_M(-2S-(p+1)\Pinf)\cong \T_{\bt^{(1,p)}}\MgNeep,
\\
\L_{-1}\oplus\L_{(0)}^{(p)}&\cong \H^1(M,T_M(-S-(p+1)\Pinf))\cong
\T_{\bt^{(0,p)}}\MgNe,
\\
\L_{(0)}^{(p)}&\cong \H^1(M,T_M(-(p+1)\Pinf)\cong
\T_{[\Sigma,\Pinf]}\Mgep.
\end{aligned}
\end{equation}
Again, for the infinite jets we obtain
\begin{equation}\label{E:modisoi}
\T_{\bt^{(1,\infty)}}\MgNeei=\lim_{p\to\infty}
\H^1(M,T_M(-2S-p\Pinf)) \cong \L_{(0)}\oplus\L_-.
\end{equation}
Let us note the dimension formula
\begin{equation}\label{E:dimMn2}
\dim_{b^{(1,p)}}(\MgNeep)=
 \begin{cases}      3g-2+2N+p,& g\ge 1\\
                          \max{(0,N-2)}+N+p,&g=0\ .
 \end{cases}
\end{equation}
For $N\ge 2$ the first expression is valid for any genus.

\medskip
Let $\bt^{(1,p)}\in\MgNeep$ be a moduli point. Let
$\nu^{(p)}:\MgNeep\to\MgNe$ be the map forgetting the coordinates
and let $\bt=\nu^{(p)}(\bt^{(1,p)})$ be a generic point with open
neighborhood $\Wt$. For $\bt=[\mpt]$ we can construct the
Krichever-Novikov objects
\begin{equation}\label{E:locob}
\A_{\tilde b},\ \L_{\tilde b},\ \Lh_{\tilde b}, \ \gb_{\tilde b},\
\gh_{\tilde b},\ \Fl_{\tilde b},\ \text{etc.}
\end{equation}
with respect to $I=\{P_1,P_2,\ldots,P_N\}$ and $O=\{\Pinf\}$.
Recall from \cite {rSSpt} that  there are  sheaf versions of these
objects
\begin{equation}\label{E:sheaf}
\A_{\Wt},\ \L_{\Wt},\ \Lh_{\Wt},\ \gb_{\Wt},\ \gh_{\Wt},\
\Fl_{\Wt}.
\end{equation}

Similarly, we can consider the sheaf versions of the objects
introduced in \refS{kzkn}, and \refS{kzrep} of the present paper:
e.g. $\D_{\g,\Wt}$ -- the sheaf of algebras of Krichever-Novikov
differential operators of order $\le 1$, and $V_{\Wt}$ -- the
sheaf of fermion modules.

Introduce the \emph{regular subalgebras} of $\gb$ and $\L$ as
follows. Let $\A^r\subset\A$, and $\L^r\subset\L$ consist of those
elements vanishing at $P_\infty$. Introduce
$\gb^{\,r}=\g\otimes\A^r$. Observe that $\gb^{\,r}$ is a Lie
subalgebra of $\gb$ as well as of $\gh$, and $\L^r$ is also a Lie
subalgebra of $\Lh$. Denote by $\gb^{\,r}_{\Wt}$, $\L^r_{\Wt}$
etc. the corresponding sheaves.
\begin{definition}\label{D:cb}
Let $\mathcal{V}_{\Wt}$ be a sheaf of (fibrewise) representations
of $\gh_{\Wt}$.  The sheaf of {\it conformal blocks} (associated
to the representation $\mathcal{V}_{\Wt}$)
 is defined as the sheaf of coinvariants
\begin{equation}\label{E:cb}
C_{\Wt}=\mathcal{V}_{\Wt}/\gb^{\,r}_{\Wt}\mathcal{V}_{\Wt}.
\end{equation}
\end{definition}

\medskip
For $p\in\N$ or $p=\infty$ let
$\Wt^{(p)}=\left(\nu^{(p)}\right)^{-1}(\Wt)$. By pulling-back the
above sheaves over $\Wt$ via $\nu^{(p)}$ we obtain sheaves on the
open subset $\Wt^{(p)}$ of $\MgNe^{(p)}$. Starting from a
fibrewise representation $\mathcal{V}_{\Wt^{(p)}}$, the sheaf of
conformal blocks can be defined in the same way as \refE{cb} by
\begin{equation}
C_{\Wt^{(p)}}=\mathcal{V}_{\Wt^{(p)}}/\gb^{\,r}_{\Wt^{(p)}}
\mathcal{V}_{\Wt^{(p)}}.
\end{equation}
Clearly, if  $\mathcal{V}_{\Wt}$ is already defined on $\Wt$ then
$\nu^*_p(C_{\Wt})=C_{\Wt^{(p)}}$.

Of special importance is the pull-back to $\Wt^{(1)}$. Recall that
this means that we fix a set of first order jets of coordinates.
As shown in \cite[Lemma 4.3]{rSSpt} this
 fixes the Krichever-Novikov basis elements
uniquely. In particular we can choose in every of the above vector
spaces a standard basis given by these Krichever-Novikov basis
elements. By their explicit form given in \cite{rSLb} it is
obvious that they depend analytically on the moduli point. In this
way we see that over $\Wt^{(1)}$ the sheaves
\begin{equation}\label{E:sheafo}
\A_{\Wt^{(1)}},\ \L_{\Wt^{(1)}},\ \Lh_{\Wt^{(1)}},\
\gb_{\Wt^{(1)}},\ \gh_{\Wt^{(1)}},\ \Fl_{\Wt^{(1)}}.
\end{equation}
are free sheaves of trivial infinite-dimensional vector bundles
with trivializations given by the Krichever-Novikov basis
(respectively constructed out of them). Of course,  everything
remains true for $\Wt^{(p)}$ instead $\Wt^{(1)}$.

Below, we take the sheaf of fermion modules
$\mathcal{V}_{\Wt^{(1)}}$ as a representation sheaf. This will be
our model situation. In our description, the basis fermions do not
depend on moduli at all, only the Lie algebra action does via the
structure constants. We trivialize the sheaf
$\mathcal{V}_{\Wt^{(1)}}$ using these bases and obtain the
corresponding trivial vector bundle $V\times\Wt^{(1)}$, where $V$
is a standard fermion space $\Hwftm$ (\refS{kzrep}). Over a
generic point of the moduli space, the space $\gb^{\,r}V$ also
does not depend on moduli (the dependence due to the structure
constants disappears after taking the linear span of images of the
basis fermions under the action of $\gb^{\,r}$). Hence, locally,
the sheaf of conformal blocks is free and defines a vector bundle.
Moreover, we will take $\g=\gl(n)$, and $\tau$ the standard
representation of $\g$ in the $n$-dimensional vector space, see
\refS{kzrep}. This guarantees the finite-dimensionality of
conformal blocks\footnote{We are grateful to B.Feigin for this
remark}. In particular, we obtain in this case that the vector
bundle of conformal blocks is of finite
rank. Note also that in this case the  representation $\Hwftm$ is irreducible%
\footnote{The irreducibility is well-studied in the graded case
\cite[Lecture 9]{rKaRa}. The proof for the almost graded case is
similar. We will give it somewhere else.}.

Observe, that all statements and arguments of the next  section
remain true for every representation sheaf $\V$ of $\gh$ with an
arbitrary reductive $\g$ as long as the above properties of our
model situation are true.

\medskip


\subsection{Projective flat connection and generalized
   Knizhnik-Zamolodchikov equations}${ }$

For short, denote by $\mathcal V$ a sheaf
$\mathcal{V}_{\Wt^{(1)}}$ fulfilling the conditions given at the
end of the last subsection, e.g. a sheaf of fermion
representations, and let it be fixed trough-out this subsection.

Let  $\tau=(\tau_1,\ldots,\tau_m)$ denote  the moduli parameters
in a  neighbourhood of a fixed point in $\M$, i.e. local
coordinates on $\M$.  Denote by $\Sigma(\tau)$ the Riemann surface
with corresponding conformal structure.

Choose a generic point with moduli parameters $\tau_0$ in $\M$.
Let $\tau_0$ be represented by the geometric data
$(\Sigma(\tau_0),P_1,\ldots,P_N,P_\infty,
z_1^{(1)}\ldots,z_N^{(1)},z_\infty^{(1)})$, see \refE{point}. In
particular,  $\Sigma(\tau_0)$ has  a fixed conformal structure
representing the algebraic curve corresponding to the moduli
parameters  $\tau_0$. For $\tau$ lying in a small enough
neighbourhood of $\tau_0$, the conformal structure on
$\Sigma(\tau)$ can be obtained  by deforming the conformal
structure $\Sigma(\tau_0)$ in the following way.

On the Riemann surface $\Sigma(\tau_0)$, we choose a local
coordinate $w$ at $P_\infty$ which has first order jet
$w_\infty^{(1)}$. Rigorously speaking, this amounts to passing
(temporary) to the moduli space $\Wt^{(\infty)}$. Let
$U_\infty\subset{\mathbb C}$ be a unit disc with natural
coordinate $z$. After identification of $z$ with $w$ we can think
of $U_\infty$ as a subset of the coordinate
 chart at $P_\infty$. In fact, we might even
assume that the disc with radius 2 is still inside this chart. Let
$v\in\L$ be a Krichever-Novikov vector field. By restriction, it
defines a meromorphic vector field on $U_\infty$. In turn, this
vector field defines a family of local diffeomorphisms $\phi_t$ --
the corresponding local flow. For $t$ small enough they map an
annulus $U_v\subset U_\infty$ which is bounded by the unit circle
from outside to a deformed annulus in the disc with radius 2. We
take the set on the coordinate chart at $P_\infty$ which (after
the identification of $z$ with $w$) is the interior complement to
$\phi_t(U_v)$, cut it out of the curve and use $\phi_t$ to define
a glueing of $U_\infty$ to the rest of the curve along the subset
$U_v$. In this way, for every $t$ we obtain another conformal
structure. Depending on the vector field (and corresponding
diffeomorphisms) the equivalence class of the conformal structure
will change or not. But, in any case we obtain via this process
any conformal structure which is close to the given one in moduli
space (see \cite{HaMo} for further information). Moreover, the
deformation in any tangent direction in the moduli spaces
$\Wt^{(1)}$ and $\Wt^{(\infty)}$ can be realized.

In abuse of notation we assign to $\tau$ such a diffeomorphism
$d_\tau$,  $w=d_\tau(z)$, where $d_\tau$ is defined on the annulus
$U_\tau$; keeping in mind that $\tau$ uniquely defines neither
$d_\tau$ nor $U_\tau$. Here we use $w$ for the new coordinate and
$z$ the standard coordinate.

In this section, let $\A$ denote the sheaf  of Krichever-Novikov
function algebras on $\Wt^{(1)}$, as well as the corresponding
infinite-dimensional vector bundle. In this bundle, denote the
fibre over $\tau\in\Wt^{(1)}$ by $\A_\tau$. If $A$ is a section of
the bundle we write it as $A(\tau)=A_\tau$, where
$A_\tau\in\A_\tau$. Similarly, let $\A^r$, $\Lh$, $\L$, $\L^r$,
$\gh$, $\gb^{\,r}$ denote the sheaves (respectively the bundles)
of the corresponding algebras.

Let $\A^{ann}_\tau$ be the algebra of regular functions on the
annulus $U_\tau$. Embed $\A_\tau$ into $\A^{ann}_\tau$ by
restricting every function  $A_\tau$ onto the image of $U_\tau$ in
$\Sigma(\tau)$ and consider the restriction as a function of the
variable $w$. Denote the result $A_\tau(w)$. On the other hand,
restrict  $A_\tau$ onto the coordinate chart at the point
$P_\infty\in\Sigma(\tau)$, and  denote the result
$\tilde{A}_\tau(z)$. Thus, we assign to every pair
$(A_\tau,\Sigma(\tau))$, in a non-unique way, a pair
$(\tilde{A}_\tau,d_\tau)$, where
$\tilde{A}_\tau(z)=A_\tau(d_\tau(z))$ on $U_\tau$. For $0<s<1$ let
$\tilde{\A}_s$ be the algebra of regular functions in the annulus
$U_s\subseteq U_\infty$ with boundary circles of radius 1 and s
respectively. We set $\tilde A=\injlim_{s\to 1}\tilde{\A}_s$ with
respect to the natural inclusion. Clearly,
$\tilde{A}_\tau\subseteq \tilde{\A}$ in the neighborhood of
$\tau_0$ (moreover, $\tilde{A}_\tau$ gives an element of the sheaf
of germs of meromorphic functions in $z$ at the point $z=0$).
Thus, the correspondence ${A}_\tau\rightarrow\tilde{A}_\tau$ gives
an embedding of $\A$ into $\tilde{\A}$. Denote the subsheaf of
sections of $\tilde\A$ corresponding to germs of analytic
functions vanishing at $P_\infty$ by ${\tilde\A}^r$. Let
${\tilde\L}^r$ have the similar meaning with respect to $\L$.

Given a vector field $X$ on $\Wt^{(1)}$ by $\partial_X A_\tau$ we
mean a full derivative in $\tau$ of
$\tilde{A}_\tau(d_\tau^{-1}(w))$ along the vector field $X$.
 We interpret this as a differentiation
(in $\tau$) of a Krichever-Novikov function as function in the
variable $z$ taking account of the dependence of the local
coordinate $z$ on moduli. After the substitution $w=d_\tau(z)$ we
consider $\partial_X A_\tau$ as an element of the sheaf
$\tilde{\A}$.

Consider a family of local diffeomorphisms $d_\tau$ where $\tau$
runs over a disk in the space of moduli parameters. Recall that
for every $\tau$ the corresponding $d_\tau$ is nothing but the
gluing function $w=d_\tau(z)$ for $\Sigma(\tau)$. From this point
of view, the family $d_\tau$ is nothing but a family of (local)
functions depending on parameters. To stress this interpretation,
define the function $d(z,\tau)=d_\tau(z)$. Define
$\partial_Xd_\tau$ as follows:
\begin{equation}
    (\partial_Xd_\tau)(z)=\sum_i X_i(\tau)\frac{\partial
    d(z,\tau)}{\partial\tau_i}.
\end{equation}
Given a vector field $X$ on $\Wt^{(1)}$ we can assign to it a
local vector field $\rho(X):=d_\tau^{-1}\cdot\partial_Xd_\tau$ (on
$\Sigma(\tau)$) which represents the {\it Kodaira-Spencer
cohomology class} of the corresponding 1-parameter deformation
family. By adding suitable coboundary terms (which amounts to
composing the diffeomorphism $d_\tau$ with a diffeomorphism of the
disk) we obtain
\begin{equation}\label{E:kod}
  \rho(X)=d_\tau^{-1}\cdot\partial_Xd_\tau\in\L.
\end{equation}
Given the vector field $X$ we can also assign to it via the
isomorphism \refE{modisoi} a  Krichever-Novikov vector field
$e_X$. Note that we consider $X$ as vector field on
 $\Wt^{(1)}$. Hence $e_X$ is only fixed up to the addition
of elements of $\L_{-}^{(1)}$. See \refSS{ags} for the definition
of
 $\L_{-}^{(1)}$; its elements correspond to the changes of the
coordinate at $\Pif$. We call every such element $e_X$ a {\it
pull-back}  of $X$. Independently of the pull-back $e_X$ and of
the choice of $\rho(X)$ (satisfying \refE{kod}) we have
$e_X-\rho(X)\in\L^r$, hence
\begin{equation}\label{E:lift}
e_X=\rho(X)+e^r,
\end{equation}
with  $e^r\in\L^r$ ($e^r$ depends on both $e_X$ and $\rho(X)$).

\begin{proposition}\label{P:nabl} For
every section $A$ of the sheaf $\A$, and every local vector field
$X$ define $A^X$ by the relation
\begin{equation}\label{E:reg}
   \partial_XA=-(e_X).A+A^X.
\end{equation}
Then $A^X\in{\tilde\A}$, and, moreover, $A^X\in{\tilde\A}^r$ for
$A\in\A^r$.
\end{proposition}
\begin{proof} In terms of the composition of maps,
$A_\tau=\tilde{A}_\tau \circ d_\tau^{-1}$. By the chain rule,
\begin{equation}\label{E:pd}
   \partial_XA_\tau = {\tilde\partial_X \tilde{A}_\tau}
   + {{\partial_z \tilde{A}_\tau}}
   \cdot\partial_Xd_\tau^{-1}={\tilde\partial_X \tilde{A}_\tau}+
   \left(-\widetilde{\rho(X)}{\frac {\partial}{\partial z}}
 \tilde{A}_\tau\right)\cdot d_\tau^{-1},
\end{equation}
where ${\tilde\partial_X \tilde{A}_\tau}$ is a derivative along
$X$ with the assumption of independence of the local coordinate
$z$ of $\tau$, $\partial_z \tilde{A}_\tau$ is a differential of
$\tilde{A}_\tau$ (in the variable $z$), and
$\widetilde{\rho(X)}{\frac {\partial}{\partial z}}$ is the  first
order differential operator corresponding to the vector field
$\rho(X)= d_\tau^{-1}\cdot\partial_Xd_\tau$. In more detail,
$\partial_Xd_\tau^{-1}=-d_\tau^{-1}\cdot\partial_Xd_\tau \cdot
d_\tau^{-1}$, hence ${{\partial_z \tilde{A}_\tau}}
\cdot\partial_Xd_\tau^{-1}= {{\partial_z \tilde{A}_\tau}}
\cdot(-d_\tau^{-1}\cdot\partial_Xd_\tau)\cdot d_\tau^{-1}$. By
\refE{kod}, for every
 $\tau$, the $\rho(X)=d_\tau^{-1}\cdot\partial_Xd_\tau$
is an element of $\L_\tau$. Denote the corresponding first order
differential operator (in the variable $z$) by
${\widetilde{\rho(X)}\frac {\partial}{\partial z}}$. Observe that
${{\partial_z \tilde{A}_\tau}}
\cdot(d_\tau^{-1}\cdot\partial_Xd_\tau)=\widetilde{\rho(X)} {\frac
{\partial}{\partial z}}{\tilde A}$.

Making use of \refE{lift}, we replace $\widetilde{\rho(X)}{\frac
{\partial}{\partial z}}$ in \refE{pd} with the differential
operator corresponding to $e_X-e^r$ and obtain \refE{reg} with
$A^X_\tau=\tilde\partial_X{\tilde A}_\tau + e^rA_\tau\in\tilde\A$.

Assume, $A_\tau\in\A^r$ for every $\tau$, i.e.
$A_\tau(P_\infty)=0$ for $P_\infty\in\Sigma(\tau)$. We have
$e^r\in\L^r$ by definition of the pull-back in question, hence
also $e^r(P_\infty)=0$. Since $z(P_\infty)=0$,
$A_\tau(P_\infty)=0$ implies ${\tilde A}_\tau(0)=0$, for every
$\tau$, and, further on, ${\tilde\partial}_X{\tilde A}_\tau(0)=0$.
Therefore, $A^X_\tau(P_\infty)=0$, hence $A^X_\tau\in
({\tilde\A}^r)_\tau$ for every $\tau$, which completes the proof.
\end{proof}
Let $\L$ be  the sheaf of  Krichever-Novikov vector field
algebras, and $e$ be a meromorphic section of it. Turning to the
definition of $\partial_Xe$, we observe that there is also no
conventional one. In analogy with the case of functions, we could
define it using local vector fields on $U_\tau$. Another
possibility, which  we prefer here, is to define it via Leibniz
rule: i.e. for every $A\in\A$, by definition,
\[(\partial_Xe).A=\partial_X(e.A)-e.\partial_XA .
\]

\begin{proposition}\label{P:reg1}
For every section $e$ of the sheaf $\L$, and every local vector
field $X$ define $e^X$ by the relation
\begin{equation}\label{E:reg2}
   \partial_Xe=-[e_X,e]+e^X.
\end{equation}
Then $e^X\in{\tilde\L}$, and, moreover, $e^X\in\tilde{\L}^r$ for \
$e\in\L^r$.
\end{proposition}
\begin{proof}
By \refP{nabl} we have $\partial_X(e.A)=-e_X.(e.A)+(e.A)^X$ where,
$(eA)^X={\tilde\partial}_X(\widetilde{e.A})+e^r.(e.A)$. This
follows from the proof of \refP{nabl}. Similarly,
$e.\partial_XA=e.(-e_X.A+{\tilde\partial}_X\tilde{A}+e^r.A)$. All
together
\begin{equation}\label{E:vdef}
(\partial_Xe).A=-[e_X,e].A+{\tilde\partial}_X(\widetilde{e.A})
    -e.{\tilde\partial}_X\tilde{A}+[e^r,e].A.
\end{equation}
Since the objects $e$ and $A$ are global, we have
$\widetilde{e.A}=\tilde{e}.\tilde{A}$. Applying the Leibniz rule
again, we obtain ${\tilde\partial}_X(\widetilde{e.A})
    -e.{\tilde\partial}_X\tilde{A}=({\tilde\partial}_X{\tilde e}){\tilde A}$.
Since \refE{vdef} is a relation in the sheaf $\tilde\A$, we do not
distinguish between $A$ and $\tilde A$. Hence, \refE{vdef} implies
\refE{reg2} where $e^X={\tilde\partial}_X\tilde{e}+[e^r,e]\in
\tilde\L$.

If $e\in\L^r$ then $[e^r,e]\in \L^r$ since $e^r\in\L^r$ and $\L^r$
is a subalgebra. Further on, ${\tilde\partial}_X\tilde{e}(0)=0$
for the same reason as ${\tilde\partial}_X{\tilde A}(0)=0$ in the
proof of \refP{nabl}. Thus, $e^X\in\tilde\L^r$ which completes the
proof.
\end{proof}

Consider a sheaf of operators on the local sections of the sheaf
$\V$. Assume $B$ to be a local section of it. By definition,
$\partial_XB=[\partial_X,B]$, where, on the right hand side,
$\partial_X$ is a differentiation on $\V$. Here the following
cases occur: $B=u(A)$, where $u\in\g$ and $A\in\A$ and
$B=T[e]:=T(e)$, the Sugawara operator introduced by
\refE{suge}%
\footnote{In this section there will be no danger of confusion of
$T(e)$ with \refE{suga},hence,  we choose $T(e)$ to avoid
confusion with the Lie bracket.}.

For every pull-back $e_X$ of $X$ we introduce the following first
order differential operator on sections of the trivial sheaf $\V$:
\begin{equation}\label{E:conn}
 \nabla_X=\partial_X+T(e_X),
\end{equation}
where $\partial_X=\sum_iX_i(\tau)\frac{\partial}{\partial\tau_i}$.

\begin{proposition}\label{P:corr}
 $\nabla_X$ is well-defined on conformal blocks and is independent
 of a pull-back of $X$ there.
\end{proposition}

Before we prove this proposition we have to extend the operators
$u(A)$ and $T(e)$ to the case when $A$ or $e$ are local objects.
Local vector fields (functions, currents etc.) form completions of
the corresponding Krichever-Novikov objects since they have
infinite expansions over the corresponding Krichever-Novikov bases
\cite{rKNFa}.
 By \refE{modiso} and
\refP{reg1} we can restrict ourselves with expansions with only  a
finite pole order at $P_\infty$, i.e. with finitely many
components of positive Krichever-Novikov degree. Define
$\overline{V}$ as the  space of infinite formal sums in negative
direction (with respect to the fermion degree \refE{mdeg}) of
linear combinations of basis fermions. Due to the
almost-gradedness, the action of the above operators can be
extended on $\overline{V}$ with the advantage, that these
extensions do exist also for $A$ and $e$ local. Indeed, only
finitely many terms of the expansions (of the operator and of the
corresponding element in $\overline V$) contribute in the result
to the component of  a given degree.

First, we prove the following Lemma.
\begin{lemma}\label{L:normal}
For the fermion representations we have
\begin{equation}\label{E:normal}
\partial_X u(A)=u(\partial_X A)
\end{equation}
\end{lemma}
\begin{proof}
For every $u\in\g$, $A\in\A$ and every basis fermion
$\Psi=\psi_{i_1}\wedge\psi_{i_2}\ldots$ we have
$u(A)\Psi=(uA)\psi_{i_1}\wedge\psi_{i_2}\ldots+\psi_{i_1}\wedge
(uA)\psi_{i_2}\ldots +\l_1\cdot\Psi$, where in the expressions
$(uA)\psi_{i_k}$ the term with $\psi_{i_k}$ (if there is any) has
to be ignored by regularization and the last term is the counter
term coming from regularization.
 Since basis fermions do not depend on
moduli, we have
\begin{equation}
\begin{split}
     (\partial_Xu(A))\Psi
    &=\partial_X(u(A)\Psi)\\
    &=(u\partial_XA)\psi_{i_1}\wedge\psi_{i_2}\ldots+\psi_{i_1}\wedge
      (u\partial_XA)\psi_{i_2}\ldots +(\partial_X\l_1)\cdot\Psi  \\
    &=u(\partial_XA)\Psi +(\partial_X\l_1)\cdot\Psi-\l_2\cdot\Psi
  \end{split}
\end{equation}
where $\l_2\cdot\Psi$ appears due to regularization of
$u(\partial_XA)$. As long as  no regularization is necessary, the
relation \refE{normal} follows immediately. The regularization can
be easily calculated via the matrix of the operator $uA$ in the
space of sections of the holomorphic bundle involved. Let
$uA=\sum\limits_{i,j=-\infty}^\infty a_{ij}E_{ij}$ where
$\{E_{ij}|i,j\in{\mathbb Z}\}$ is the natural basis in the matrix
space.
 By the regularization procedure \cite{rKaRa}
\begin{equation*}
\l_1=\sum_{i\in N_-} a_{ii} -\sum_{i\in N_+} a_{ii},
\end{equation*}
where $N_+$ is the set of non-occupied positions, or holes of
positive degree in $\Psi$, i.e. ${N_+}={\mathbb N}
\setminus\{i_1,i_2,\cdots\}$, and
$N_-$ is the set of occupied  positions of degree $\le 0$.%
\footnote{Other descriptions are possible, but they are
equivalent.} Similarly, $\l_2=\sum_{i\in N_-} \partial_Xa_{ii}
-\sum_{i\in N_+} \partial_Xa_{ii}$.
 Therefore,
$\partial_X\l_1-\l_2=0$ and the claim is true also in this case.
\end{proof}

For more general  representations, we take relation \refE{normal}
as an additional requirement.

\begin{proof}[Proof of  \refP{corr}]
Using \refL{normal} and \refL{fund} we find
\[
\begin{split}
[\nabla_X,u(A)]&=[\partial_X+T(e_X),u(A)] =
[\partial_X,u(A)]+[T(e_X),u(A)] \\
&= u(\partial_XA)+u(e_X.A).
\end{split}
\]
Assume, $A\in{\A^r}$. Then, by \refP{nabl}, we have
\[ u(\partial_XA)=-u(e_X.A)+u(A^X),
\]
where $A^X\in\A^r$. Hence, $[\nabla_X,u(A)]=u(A^X)$ and
\[ [\nabla_X,u({\A^r)}]\subseteq
u({\A^r)}.
\]
Hence $\gb^r\V$ is a $\nabla_X$-invariant subspace and $\nabla_X$
is well-defined on $\V/\gb^r\V$.
\end{proof}
\begin{lemma}\label{L:norm1}
For every $X\in T\M$ we have
  $$\partial_XT(e)=T(\partial_Xe)+\l\cdot id,
  $$
where $\l=\l(X,e)\in{\C}$.
\end{lemma}
\begin{proof}
By the fundamental relation (\refL{fund}) for every $e\in\L$,
$u\in\g$, $A\in\A$
\begin{equation}\label{E:bas}
[T(e),u(A)]=u(e.A).
\end{equation}
Take the derivative on both sides of the relation \refE{bas} along
a local vector field $X\in \T\Wt^{(1)}$. By \refL{normal}  we
obtain
\begin{equation}\label{E:dr}
[\partial_XT(e),u(A)]+[T(e),\partial_Xu(A)]=
    u((\partial_Xe).A)+u(e.(\partial_XA)).
  \end{equation}
  Again, by \refE{bas} and \refL{normal} the second terms on
both sides of \refE{dr} are equal. Therefore,
  $$[\partial_XT(e),u(A)]=u((\partial_Xe).A).
  $$
Applying \refE{bas} once more, we replace the right hand side of
the latter relation by $[T(\partial_Xe),u(A)]$ (see the remark
below). Therefore,
 \begin{equation}\label{E:com}
[\partial_XT(e)-T(\partial_Xe),u(A)]=0
\end{equation}
for every  $X\in \T\Wt^{(1)}$, $u\in\g$, $A\in\A$.

By standard arguments of the theory of highest weight
representations (either by the irreducibility of  the
representation or by uniqueness of the vacuum vector), the
commutation relations \refE{com} immediately imply the lemma.
\end{proof}

\begin{remark}
The $\partial_Xe$ is a local vector field on a deformed annulus
(an element of the sheaf $\tilde\L$ to be more precise). Hence, we
need the relation \refE{bas} for local vector fields to prove
\refE{com}. Due to the definition given after the formulation of
\refP{corr}, the representations $T(e)$ and $u(A)$ are well
defined also  on local vector fields and functions, respectively,
with preserving the relation \refE{bas}. Indeed, for an $A\in\A$,
a homogeneous $v\in V$ and an arbitrary $n$ there exists a partial
sum $\tilde{e}$ of the expansion for $\partial_Xe$ such that
$(u((\partial_Xe).A)v)_n=(u({\tilde e}.A)v)_n$ and
$([T(\partial_Xe),u(A)]v)_n=([T({\tilde e}),u(A)]v)_n$ where
$(\cdot)_n$ denotes the projection onto the component of degree
$n$ in $\overline V$. By \refE{bas} the right hand sides of the
last two relations are equal, hence their left hand sides also are
equal. This implies the relation \refE{bas} with $\partial_Xe$
instead $e$.
\end{remark}

For the remainder of this section, our goal is to prove the
projective flatness of \refE{conn} and to introduce the
corresponding analog of the Knizhnik-Zamolodchikov equations.

\begin{lemma}\label{L:ue} For every pull-backs $e_X$, $e_Y$
of  local vector fields $X$, $Y$ to $\L$, there exist a pull-back
$e_{[X,Y]}$ of $[X,Y]$ such that
\[ e_{[X,Y]}=[e_X, e_Y]+\partial_Xe_Y-\partial_Ye_X.
\]
\end{lemma}
\begin{proof} By \cite[Lemma 1.3.8]{rTUY}
\begin{equation}\label{E:rho}
\rho([X,Y])=[\rho(X),\rho(Y)]+\partial_X\rho(Y)-\partial_Y\rho(X).
\end{equation}
Observe that for every pull-backs $e_X$, $e_Y$ we have
$e_X=\rho(X)+e_1^r$, $e_Y=\rho(Y)+e_2^r$, where
$e_{1}^,re_{2}^r\in\L^r$, see \refE{lift}. In terms of $e_X$,
$e_Y$  \refE{rho} reads as
\begin{equation}\label{E:rho1}
\rho([X,Y])=[e_X,e_Y]+\partial_Xe_Y-\partial_Ye_X+e_3^r
\end{equation}
where $e_3^r= [e_1^r,e_2^r]-(\partial_Xe_2^r+[e_X,e_2^r])+
(\partial_Ye_1^r+[e_Y,e_1^r])$.

Since $\L^r$ is a Lie subalgebra, we have $[e_1^r,e_2^r]\in\L^r$.
By \refP{reg1}, the elements $\partial_Xe_2^r+[e_X,e_2^r]$ and
$\partial_Ye_1^r+[e_Y,e_1^r]$ are also regular, hence
$e_3^r\in\L^r$. Thus, $e_{[X,Y]}=\rho([X,Y])-e_3^r$ is a pull-back
of $[X,Y]$, and the lemma is proven.
\end{proof}

\begin{theorem}\label{T:flat}
$\nabla_X$ is a projective flat connection on the vector bundle of
conformal blocks.
\end{theorem}\begin{proof}
\begin{equation}\label{E:fl1}
\begin{split}
     [\nabla_X,\nabla_Y]
    &=[\partial_X+T(e_X),\partial_Y+T(e_Y)]\\
    &=[\partial_X,\partial_Y]+[\partial_X,T(e_Y)]-
      [\partial_Y,T(e_X)]+[T(e_X),T(e_Y)].
    \end{split}
  \end{equation}
  Since $T$ is a projective representation of $\L$ and due to the
relations $[\partial_X,T(e_Y)]=\partial_XT(e_Y)$,
$[\partial_Y,T(e_X)]=\partial_YT(e_X)$, we can rewrite \refE{fl1}
in the following form:
\begin{equation}\label{E:fl2}
 [\nabla_X,\nabla_Y]
     =\partial_{[X,Y]}+T(\partial_Xe_Y- \partial_Ye_X+[e_X,e_Y])
       +\l\cdot id.
     \end{equation}
Here we used also \refL{norm1}. By \refL{ue}, this reads as
\begin{equation}\label{E:fl3}
\begin{split}
     [\nabla_X,\nabla_Y]
     &=\partial_{[X,Y]}+T(e_{[X,Y]})+\l\cdot id\\
     &=\nabla_{[X,Y]}+\l\cdot id.
   \end{split}
 \end{equation}
\end{proof}
We consider the following equations for horizontal sections of the
connection $\nabla_X$ as a generalization of
Knizhnik-Zamolodchikov equations:
\begin{equation}\label{E:KZ}
       \nabla_X\Psi=0,\quad X\in H^0(\mathcal{U}, \T\MegN)
\end{equation}
where $\Psi$ is a section of the sheaf of conformal blocks. These
equations are proposed in \cite{rSSpt}, also $\nabla_X$ was
explicitly calculated there for $g=0$ and $g=1$.


\end{document}